\documentclass[10pt,numreferences]{kluwer}
\usepackage{amssymb,epic,eepic}

\font\fr=eufm10  scaled \magstep 1
\def\cal#1{\mathcal{#1}}

\def\l{\lambda}
\def\a{\alpha}

\def\t{\tau}
\def\s{\sigma}
\def\w{\omega}

\def\G{\Gamma}
\def\g{\gamma}
\def\d{\delta}
\def\e{\epsilon}
\def\T{{\mathbf T}}
\def\an{\Lambda} 
\def\cinfty#1{C^{\scriptscriptstyle\infty}(#1)}
\def\vectorfields#1{{\cal X}(#1)}
\def\vvectorfields#1{{\cal V}(#1)}
\def\voneforms#1{{\cal V}^*(#1)}

\def\fpd#1#2{\frac{\partial #1}{\partial #2}}
\def\R{\dR}
\def\N{\dN}

\def\del{\nabla}
\def\J{{\cal J}}
\def\ovl#1{\overline{#1}}

\def\im{\mbox{Im }}

\def\oneforms#1{{\cal X}^*(#1)}

\def\ro{\rho}

\def\een{\mbox{id}}
\def\sign{\mathrm{sgn}}

\newtheorem{thm}{Theorem}
\newtheorem{rem}[thm]{{\normalfont{\it Remark}}}
\newtheorem{prop}[thm]{Proposition}
\newtheorem{defn}[thm]{Definition}
\newtheorem{cor}[thm]{Corallary}
\newtheorem{lem}[thm]{Lemma}

\def\etal{{et al.}}
\begin{document}
\begin{article}
\begin{opening}
\title{Geometric aspects of the Maximum Principle and lifts over a bundle map}

\author{B. \surname{Langerock}\email{bavo.langerock@rug.ac.be}}
\institute{Department of Mathematical Physics and Astronomy, Ghent
University, Krijgslaan 281 S9,B-9000 Gent Belgium}

\runningtitle{Geometric aspects of the Maximum Principle}

\begin{abstract}
A coordinate-free proof of the Maximum Principle is provided in
the specific case of an optimal control problem with fixed time.
Our treatment heavily relies on a special notion of variation of
curves that consist of a concatenation of integral curves of
time-dependent vector fields with unit time component, and on the
use of a concept of lift over a bundle map. We further derive
necessary and sufficient conditions for the existence of so-called
abnormal extremals.
\end{abstract}

\classification{AMS}{49Kxx, 53Cxx.}
\end{opening}

\section{Introduction and preliminary definitions}
\label{sectieintro} The results presented in this paper find their
origin in some recent work on sub-Riemannian geometry \cite{sub},
and are also strongly inspired by some ideas developed in the book
by L.S. Pontryagin \etal \ \cite{Pont}. The main purpose is to
provide a comprehensive and coordinate-free proof of the Maximum
Principle and, at the same time, to present a version of this
principle that may be readily accessible to researchers studying
the variational approach to dynamical systems subjected to
nonholonomic constraints, also called Vakonomic dynamics.
Applications of our results can be found, for instance, in
sub-Riemannian geometry, where the problem of characterizing
length-minimizing curves (see \cite{sub} and references therein)
can be solved by means of the Maximum Principle. Also the
construction of a Lagrangian and Hamiltonian dynamics on
Lie-algebroids (see, for instance, \cite{marti,martiacta,wein})
could be tackled using the formalism described in the present
paper. This will discussed in a forthcoming paper.

For the present paper, we restrict ourselves to control problems
satisfying strong smoothness conditions and we only consider
optimal control problems with fixed time. The extension of our
results to more general optimal control problems is currently
under investigation.

We now first give some preliminary definitions and fix some
notations. All manifolds considered in this paper are real, finite
dimensional smooth manifolds without boundary, and by smooth we
will always mean of class $C^{\scriptscriptstyle\infty}$. The set
of (real valued) smooth functions on a manifold $B$ will be
denoted by $\cinfty{B}$, the set of smooth vector fields by
$\vectorfields{B}$ and the set of smooth one-forms by
$\oneforms{B}$. The set of all smooth (local or global) sections
of an arbitrary fibre bundle $\tau: E \rightarrow B$ will be
denoted by $\G(\tau)$. A family ${\cal D}$ of vector fields on a
manifold $B$ is said to be \emph{everywhere defined} if, given any
point $x \in B$, there exists an element $X \in {\cal D}$ such
that $x$ is contained in the domain of $X$.

We now recall the concept of piecewise curve as introduced in
\cite{sub}. First of all, by a \emph{curve} in an arbitrary
manifold $B$ we shall always mean a smooth mapping $c: I
\rightarrow B$, with $I \subset \R$ a closed interval, and such
that $c$ admits a smooth extension to an open interval containing
$I$. A mapping $c:[a,b] \rightarrow B$ will be called a
\emph{piecewise curve} in $B$ if there exists a finite subdivision
$a_0:=a < a_1 < \ldots < a_{\ell-1} < a_{\ell}: =b$ such that the
following conditions are fulfilled:
\begin{enumerate} \item $c$ is left continuous at each point $a_i$
for $i =1, \ldots, \ell$, i.e. $\lim_{t \to a_i^{-}}c(t)$ exists and equals
$c(a_i)$;
\item $\lim_{t \to a_i^{+}}c(t)$ is defined for all $i = 1,
\ldots, \ell$ and $\lim_{t \to a_0^{+}}c(t)=c(a_0)$ (i.e.\ $c$ is right
continuous at $a_0=a$);
\item for each $i = 1,\ldots,\ell$, the mapping $c^i: [a_{i-1},a_{i}]
\rightarrow B$, defined by $c^i(t) = c(t)$ for $t \in ]a_{i-1},a_{i}]$ and
$c^i(a_{i-1}) = \lim_{t \to a_{i-1}^{+}}c(t)$, is smooth (i.e. is a curve in
$B$).
\end{enumerate} We will also say that the piecewise curve $c$ is
``induced by the smooth curves $c^i$". A piecewise curve which is
continuous everywhere will simply be called a \emph{continuous
piecewise curve} and it corresponds to what is usually called a
`piecewise smooth curve' in the literature. For example, consider
two smooth curves $\g^i:[a_{i-1},a_i] \to B$ with $i=1,2$ such
that $\g^1(a_1)=\g^2(a_1)$. According to the above definition, the
curve $\g:[a_0,a_2] \to B$ defined by $\g(t) = \g^i(t)$ if $t \in
]a_{i-1},a_i]$ and $\g(a_0) = \g^1(a_0)$, is a continuous
piecewise curve induced by $\g^1,\g^2$. On the other hand, the
piecewise curve $\dot \g$, induced by $\dot \g^1,\dot \g^2$,
provides an example of a piecewise curve which, in general, need
not be continuous.

In this paper we will also encounter the notion of piecewise
section of a bundle fibred over the real line, say $\pi: B \to
\R$, the definition of which is similar to that definition of a
piecewise curve. A smooth section $\s \in \G(\pi)$, defined on a
closed interval $I = [a,b]$, is always assumed to be the
restriction of a smooth section of $\pi$ defined on an open
interval containing $I$. Clearly, any section of $\pi$ determines
a curve in $B$. On the other hand, if $\g: I \to B$ is a curve in
$B$, then it will determine a section of $\pi$ iff $\pi(\g(t)) =
t$ for all $t \in I$. We say that $\s: I = [a,b] \to B$ is a {\em
piecewise section} of $\pi$ if $\s$ is a piecewise curve in $B$
and, in addition, $\pi(\s(t)) =t$ for all $t \in I$. Let
$\s^i:[a_{i-1},a_i] \to B$, with $i=1,\ldots,\ell$ and $a_0=a <a_1
< \ldots < a_\ell=b$, represent a finite number of curves that
induce such a piecewise section $\s$. Then, the curves $\s^i$
necessarily satisfy $\pi(\s^i(t))=t$, which implies that they are
smooth (local) sections of $\pi$. We then say that the smooth
sections $\s^i$ induce the piecewise section $\s$. A {\em
continuous piecewise section} $\s$ is a piecewise section $\s: I
\to B$ such that, in addition, $\s$ is a continuous mapping.

\section{A geometric framework for control theory}
\label{sectiegeometricsetting} We can now proceed towards the
construction of a differential geometric setting for certain
control problems. It should be emphasized that, although our
formulation is not the most general one, if only for the rather
strong smoothness conditions we impose, it occurs to us that there
is a sufficiently large and relevant class of control problems
that fit within the framework described below (see for instance
\cite{sussmann4} for a different approach).

\begin{defn}
A {\rm geometric control structure} is a triple $(\t,\nu,\ro)$
consisting of (i) a fibre bundle $\t: M \to \R$ over the real
line, where $M$ is called the event space, (ii) a fibre bundle
$\nu : U \to M$, called the control space, and (iii) a bundle
morphism $\ro: U \to J^1\t$ over the identity on $M$, such that
$\t_{1,0} \circ \ro = \nu$.
\end{defn}
In the above, $J^1\t$ is the first jet bundle of $\t: M \to \R$,
with projections $\t_1: J^1\t \to \R$ and $\t_{1,0}: J^1\t \to M$.
The typical fibre of $M$ plays the role of configuration space and
will be denoted by $Q$. It follows from the definition that we
have the following commutative diagram:

\begin{center}\setlength{\unitlength}{.8cm}
\begin{picture}(4.73,5)
\put(0.996,0.297){\parbox[t]{0.487\unitlength}{$\R$}}
\put(1.17,2.14){\vector(0,-1){1.40}}
\put(0.975,2.44){\parbox[t]{0.551\unitlength}{$M$}}
\put(1.19,4.22){\vector(0,-1){1.34}}
\put(1.04,4.51){\parbox[t]{0.445\unitlength}{$U$}}
\put(1.65,4.68){\vector(1,0){1.72}}
\put(3.69,4.66){\parbox[t]{0.911\unitlength}{$J$$^{1}$$\tau{}$}}
\put(3.57,4.25){\vector(-4,-3){1.91}}
\put(2.42,5.02){\parbox[t]{0.381\unitlength}{$\rho{}$}}
\put(0.699,1.48){\parbox[t]{0.318\unitlength}{$\tau{}$}}
\put(0.699,3.67){\parbox[t]{0.403\unitlength}{$\nu{}$}}
\put(3.12,3.43){{$\tau{}$$_{1,0}$}}
\end{picture}\end{center}

Let $u$ denote a (local) section of $\t \circ \nu$, i.e. $u
:I\subseteq \R \to U$ with $\t(\nu(u(t)))=t$. With $u$ we can
associate a section $c$ of $\t$, called the {\em base section of
$u$} and defined by $c= \nu \circ u$.
\begin{defn} A smooth section $u \in \G(\t \circ \nu)$ is said to be a
{\rm smooth control} if $\ro\circ u = j^1 c$, where $c$ denotes
the base section of $u$ and $j^1 c$ its first jet extension. A
smooth section $c \in \G(\t)$ is called a smooth controlled
section if $c$ is the base section of a smooth control $u$.
\end{defn}
Let $(t,x^i,u^a)$ denote an adapted coordinate system on $U$
(i.e.\ adapted to both fibrations $\t$ and $\nu$). The condition
for $u \in \G(\t \circ \nu)$ to be a smooth control is expressed
in coordinates as follows: putting $u(t) = (t,x^j(t),u^a(t))$ we
must have that $\ro^i(t,x^j(t),u^a(t)) = \dot x^i(t)$ for all $t$.
Note that these equations are in agreement with the definition of
a control as given in \cite[p 56]{Pont}, where $M=\R \times \R^n$
and $U$ is an (open) subset of $M\times \R^k$.

\begin{defn} A {\rm control} $u: I=[a,b] \to U$ is a piecewise section
of $\t \circ \nu$ such that $u$ is induced by a finite number of
smooth controls and, in addition, its projection $\nu \circ u$ is
a continuous piecewise section of $\t$. A continuous piecewise
section $c: I \to M$ of $\t$ will be called be a {\rm controlled
section} if it is the base section of a control.\end{defn}

In the following, we shall show that one can associate with any
section of $\nu$ a vector field on $M$. These vector fields will
generate controls in the sense that (segments of) their integral
curves can be regarded as controlled sections of $\t$. Moreover,
we shall see that also the converse holds: each controlled section
appears to consist of a concatenation of integral curves of such
vector fields. First, we shall specify what we precisely mean by a
``concatenation of integral curves" of vector fields.

Let $B$ denote an arbitrary manifold and consider a finite {\em
ordered} set of, say, $\ell$ vector fields on $B$:
$(X_{\ell},\ldots,X_1)$, which need not all be different. Let
$\{\phi^i_s\}$ denote the flow of $X_i$. The {\em composite flow}
$\Phi$ induced by $(X_\ell,\ldots,X_1)$ is then defined as the
mapping
\[ \Phi: V \subset \R^\ell \times B \to B: ((t_\ell ,\ldots, t_1),x) \mapsto
\phi^\ell_{t_\ell} \circ \ldots \circ \phi^1_{t_1}(x)\,,
\]
whose domain is a subset $V$ of $\R^\ell \times B$. For brevity we
shall write $\Phi_T(x)$ for $\Phi ((t_\ell,\ldots,t_1),x)$, where
$T:=(t_\ell,\ldots,t_1)$. We shall sometimes refer to $T$ as the
{\em composite flow parameter}. Assume that $(t_1,x) \in
\mbox{Dom} (\phi^1)$ and that $(t_{i+1},(\phi^i_{t_i}\circ \ldots
\circ\phi^1_{t_1})(x)) \in \mbox{Dom} (\phi^{i+1})$ for
$i=1,\ldots,\ell-1$, then $((t_\ell, \ldots,t_1),x) \in
\mbox{Dom}(\Phi)$. It can be proven that $\mbox{Dom}(\Phi)(= V)$
is an open set (which might be empty) and that for each $x \in B$,
$T\mapsto \Phi_T(x)$ is a smooth mapping defined on an open
neighborhood of $0 \in \R^{\ell}$. If we fix a value $T \in
\R^{\ell}$ of the composite flow parameter, then $\Phi_T: B \to B$
determines a diffeomorphism defined on an open subset of $B$. We
refer to \cite{Lib} (Appendix 3) for further details on composite
flows.

Fixing again some $T =(t_{\ell},\ldots,t_1) \in pr_1(V) \subset
\R^\ell$ (with $pr_1$ the projection of $V$ onto $\R^{\ell}$), we
can associate with any $x \in \mbox{Dom}(\Phi_T)$ and with
arbitrary $a_0 \in \R$, a continuous piecewise curve $\g:
[a_0,a_0+|t_1|+\ldots +|t_\ell|]\to B$ as follows: putting $a_i =
a_0+ \sum_{j=1}^i|t_j|$ and ${\sign(t_i)} := 0, +1, -1 $ depending
on whether $t_i = 0, t_i > 0, t_i < 0$, respectively, let
\[
\g(t) = \left\{\begin{array}{lll} \phi^1_{\sign(t_1)(t-a_0)}(x) &
\mbox{ for } & t \in [a_0,a_1]\\
\phi^2_{\sign(t_2)(t-a_1)}(\phi^1_{t_1}(x))&\mbox{ for } &
t\in\ ]a_1,a_2]\\
\ldots& & \\
\phi^\ell_{\sign(t_\ell)(t-a_{\ell-1})}(\ldots
\phi^2_{t_2}(\phi^1_{t_1}(x)))&\mbox{ for } & t\in \
]a_{\ell-1},a_\ell],\end{array}\right.
\]
For $t \in ]a_{i-1},a_i[$ we then have $\dot \g(t) = \sign(t_i)
X^i(\g(t))$ and, hence, the restriction of $\g$ to $]a_{i-1},a_i[$
is an integral curve of $X_i$, resp.\ $-X_i$, for $t_i > 0$,
resp.\ $t_i <0$. Note that $\g(a_\ell) = \Phi_T(x)$, i.e.\ the
endpoint of $\g$ coincides with the image of $x$ under the
composite flow map $\Phi_T$. If all $t_i \ge 0 $, then we say that
$\g$ is a {\em concatenation of integral curves} through $x$
associated with $\Phi$ (or, with the ordered set
$(X_\ell,\ldots,X_1)$) and corresponding to the value $T$ of the
composite flow parameter. Indeed, we than have $\dot \g (t) =
X_i(\g(t))$ for any $t \in ]a_{i-1},a_i]$.

Let us now return to the geometric control structure
$(\t,\nu,\ro)$ and recall the definition of the total time
derivative operator on the first jet bundle $J^1\t$, denoted by
$\T : J^1\t \to TM$. This is the vector field along the projection
$\t_{1,0}$ defined by
\[\T(j^1_tc) = T c\left(\left. \fpd{}{t}\right|_t\right),\]where $c \in
\G(\t)$. Note that $\t_M (\T(j^1_tc)) = c(t) = \t_{1,0}(j^1_tc)$.
Let $\s$ be a section of $\nu$, then $\ro \circ \s$ is a section
of $\t_{1,0}$ and composing it with the total time derivative, we
obtain a mapping $\T \circ \ro \circ \s: M \to TM$, which is a
smooth section of $\t_M$. The vector field $\T \circ \ro \circ \s$
is projectable with respect to $\t$, and its projection on $\R$ is
given by $\fpd{}{t}$, i.e. $T \t \circ \T \circ \ro \circ \s =
\fpd{}{t} \circ \t$. This implies that, if $\{\phi_s\}$ denotes
the flow of $\T \circ \ro \circ \s$ and $\{\l_s\}$ the flow of
$\fpd{}{t}$ on $\R$ (i.e. $\l_s(t) = t + s$), then the equality
$(\t \circ \phi_s)(m) = \l_s(\t(m))$ holds for any $m \in M$ and
for all $s$ in a neighborhood of $0$ such that $\phi_s(m)$ is
defined.

For a given $\s \in \G(\nu)$, let $\{\phi_s\}$ again denote the
flow of the vector field $\T \circ \ro \circ \s$ on $M$. Assume
that $m \in \mbox{Dom}(\phi_{\epsilon})$, for some fixed $\epsilon
> 0$ and let $\t(m) = t_0$. Consider then the curve $\g(t) =
\phi_{t-t_0}(m)$ in $M$, defined on $I=[t_0,t_0+\e]$. From the
above we know that $\t(\g(t)) = \l_{t-t_0}(t_0) = t$, implying
that $\g: I \to M$ can be regarded as a section of $\t$. Moreover,
$\g$ is the base section of the section $u: I \to U$ defined by
$u(t) = \s(\g(t))$. From the definition of $u$ it easily follows
that $\T(\ro(u(t))) = \dot \g(t) = \T(j^1_t\g)$, which is
equivalent to $\ro(u(t))=j^1_t\g$ and, hence, $u$ is a smooth
control. We may therefore conclude that (up to a
reparameterization) any integral curve of $X=\T \circ \ro \circ
\s$, with $\s \in \G(\nu)$, determines a controlled section.
Indeed, if $\g:[a,b]\to M$ is such an integral curve, with
$m=\g(a)$, then the curve $\g': [\t(m),\t(m) +b-a] \to M\,,t
\mapsto \g(t-\t(m)+a)$ is a reparametrization of $\g$,
representing a controlled section of $\t$. From now on, it will
always be tacitly assumed that the integral curves $\g$ of a
vector field of the form $\T \circ \ro \circ \s$ will be
parameterized in this way.

We now introduce the following everywhere defined family of vector
fields on $M$:
\begin{eqnarray}\label{distr}
{\cal D} = \{\T \circ \ro \circ \s \ | \ \forall \s \in \G(\nu)\}.
\end{eqnarray}
Take some arbitrary sections $\s_i \in \G(\nu), i=1,\ldots,\ell$
and put $X_i: = \T \circ \ro \circ \s_i \in {\cal D}$. Then, any
concatenation $\g:[a_0,a_\ell]\to M$ of integral curves associated
to the ordered set $(X_\ell, \ldots,X_1)$, corresponding to a
value parameter $T \in \R^\ell_+$ of the composite flow parameter
(where $\R^\ell_+ = \{ (t_\ell,\ldots,t_1) \ | \ t_i \ge 0\}$) and
such that $a_{\ell} = a_0 + t_1 + \ldots + t_{\ell}$, determines a
piecewise controlled section if $\t(\g(a_0))=a_0$. Indeed, it is
an easy exercise to see that the piecewise section $u$ induced by
$u_i(t) = \s_i(\g(t))$, for $t \in [a_{i-1},a_i]$, controls $\g$
(we are using here the notations of Section \ref{sectieintro}). In
the following we prove that the converse also holds, i.e. the base
section of any control can be regarded as a concatenation of
integral curves of vector fields belonging to ${\cal D}$. We only
prove the result for smooth controls; the proof for the more
general case then easily follows.

Let $u: I \to U$ be a smooth control with base section $c:= \nu
\circ u$. First, assume that the image $u(I)$ is contained in the
domain of an adapted coordinate chart $V$ of $U$ with coordinates
$(t,x^i,u^a)$. Consider a smooth extension $\tilde u$ of $u$,
defined on an open interval $\tilde I$ containing $I$, i.e.
$\tilde u: \tilde I \to U$ is a local section of $\t \circ \nu$
with $\tilde u(t) = u(t)$ for all $t \in I$. Upon reducing $\tilde
I$ if necessary, we may always assume that $\tilde{u}(\tilde{I})
\subset V$, and in terms of the adapted coordinates on $V$ we can
then write $\tilde u(t) = (t,x^i(t), \tilde u^a(t))$. We can now
define a local section $\s$ of $\nu$ on the open subset $V'=\nu(V)
\cap \t^{-1}(\tilde I)$ of $M$ as follows: $\s(t,x^i) =
(t,x^i,\tilde u^a(t))$, $\forall (t,x^i) \in V'$. The map $\ro
\circ \s$ determines a section of $\t_{1,0}$ satisfying $\ro \circ
\s(c(t)) = j^1c(t)$ for any $t \in I$. This implies that $c$ is an
integral curve of $\T \circ \ro \circ \s$. In case the image set
$u(I)$ is not fully contained in an adapted coordinate chart, we
can always cover the compact set $u(I)$ with a finite number of
adapted coordinate charts and choose a subdivision of $I$ such
that the image of each subinterval is entirely contained in one of
these coordinate charts. The construction above can then be
carried out for the restriction of $u$ to each of these
subintervals, and it readily follows that the base section $c$ is
a concatenation of integral curves of vector fields in ${\cal D}$.
As mentioned above, the extension of this proof to the case of
general controls is straightforward. Summarizing, we have shown
that the following property holds.

\begin{prop} \label{kernprop}
A continuous piecewise section $c: I \to M$ is a controlled
section iff $c$ is a concatenation of integral curves of vector
fields in ${\cal D}$.
\end{prop}
With the family of vector fields ${\cal D}$ on $M$ we can
associate a `quasi-order relation' $R$ on $M$ (i.e. a reflexive
and transitive relation) as follows: $R$ is the subset of $M
\times M$ defined by $(m,n)\in R$ if there exists a control
$u:[a,b]\to U$ such that $\nu(u(a))=m$ and $\nu(u(b))=n$ (we will
say that `the control $u$ takes $m$ to $n$'). For brevity we shall
also denote $(m,n)\in R$ by $m\to n$, and if we want to indicate
the control $u$ explicitly, we will write $m\stackrel{u}{\to}n$.
From Proposition \ref{kernprop}  it follows that $m \to n$ iff
there exists a composite flow $\Phi$ associated with an ordered
set $(X_\ell,\ldots , X_1)$, with $X_i \in {\cal D}$, such that
$n=\Phi_T(m)$ for some $T \in \R^\ell_+$. For any $m \in M$, the
subset $R_m \subset M$, defined by
\[
R_m = \{n \in M \ | \ m\to n\}\,,
\]
is called the {\em set of reachable points from $m$}.

In the next section we will first show that a quasi-order relation
can be associated to any everywhere defined family of vector
fields on an arbitrary manifold, and that the notion of `set of
reachable points' can be introduced in this more general setting.
We will then investigate some properties of a set of reachable
points that will play an important role in the further treatment.

\section{Some properties of the set of reachable points}
\label{sectie3}

Given an everywhere defined family of vector fields ${\cal D}$ on
an arbitrary manifold $B$ one can define a quasi-order relation
$R$ on $B$ as follows: for $x, y \in B$ we put $(x,y) \in R$ if
there exists a composite flow $\Phi$, associated with an ordered
set $(X_\ell,\ldots, X_1)$, $X_i \in {\cal D}$, such that
$\Phi_T(x) = y$ for some $T \in \R^\ell_+$. We then also write $x
\stackrel{(\Phi,T)}{\longrightarrow} y$ (or simply $x \to y$). As
described in the previous section, the concatenation of integral
curves through $x$, determined by $\Phi$ and $T$, is a continuous
piecewise curve $\g$ such that $\dot \g(t) \in {\cal D}$ for all
$t$ where the derivative exists. As in the previous section, we
can then define the set of reachable points from $x \in B$ as the
subset $R_x = \{y \in B \ | \ x \to y \}$. Note that $R_x \neq
\emptyset$ for all $x \in B$, since ${\cal D}$ is assumed to be
everywhere defined.

Let $D$ denote the smallest generalized integrable distribution,
generated by ${\cal D}$ (in the sense of H.J. Sussmann, see e.g.
\cite{Lib}) and let us denote the leaf of $D$ through a given
point $x \in B$ by $L_x$. Recall that $D_x$ is defined as the
space spanned by all tangent vectors of the form $T\Phi_{T}
(Y((\Phi_T)^{-1}(x)))$, for $Y\in {\cal D}$, $\Phi$ a composite
flow associated with a (finite) ordered set of vector fields
belonging to ${\cal D}$, and $T\in \R^\ell$ such that $x \in
\hbox{Im}(\Phi_T)$. Then it is a simple exercise to see that $R_x
\subset L_x$ for any $x \in B$. If ${\cal D} = -{\cal D}$, then
the relation $R$ is symmetric. Indeed, if $\Phi_T(x) = y$, with
$\Phi$ the composite flow determined by $(X_\ell,\ldots, X_1)$ and
$T=(t_\ell,\ldots,t_1) \in \R^\ell_+$, then $(\Phi_T)^{-1}(y) = x$
and an elementary computation shows that $(\Phi_T)^{-1} =
\Psi_{T^*}$, where $\Psi$ is the composite flow corresponding to
the ordered set $(-X_1, \ldots , -X_\ell)$ (where, by assumption,
$-X_i \in {\cal D}$) and $T^* = (t_1, \ldots , t_\ell)$, i.e.\ we
also have $y\to x$. In this case $R$ determines an equivalence
relation for which the equivalence classes are precisely the leafs
of the foliation of the smallest integrable distribution $D$
generated by $\cal D$, i.e. $R_x =L_x$ for any $x \in B$.

\begin{rem}\label{control} \textnormal{It should be emphasized here
that the everywhere defined family of vector fields (\ref{distr})
associated to a control structure, can never be invariant under
multiplication by $-1$ since, by construction, each vector field
belonging to this ${\cal D}$ is of the form $\T \circ \rho \circ
\s$ for some $\s \in \G(\nu)$ and, therefore, projects onto the
fixed vector field $\fpd{}{t}$ on $\R$. Moreover, the relation $m
\to n$ is an order relation (i.e. transitive, reflexive and not
symmetric) since, if $m \to n$ then $\t(m) \le \t(n)$
holds.}\end{rem}

We will now investigate the local structure of the set of
reachable points $R_x$ for a given $x \in B$. For that purpose we
will introduce a special class of variations of a concatenation of
integral curves of vector fields in $\cal D$, connecting $x$ with
some $y \in R_x$, such that these variations will lead us from $x$
to points in a neighborhood of $y$ that also belong to $R_x$. The
following description is merely intended to give a general
intuitive idea of the kind of variation we have in mind. We will
be more specific later on.

Consider the composite flow $\Phi$ corresponding to an ordered set
of, say, $\ell$ vector fields in ${\cal D}$, and let $T \in
\R^\ell_+$ be such that $\Phi_T(x)=y$. Let $\g: [a,b] \rightarrow
B$ be the concatenation of integral curves induced by $\Phi$ and
$T$, as constructed in the previous section, with $\g(a)=x$ and
$\g(b)=y$. Roughly speaking, a {\em variation of $\g$} consists of
a $1$-parameter family of continuous piecewise curves $\g_{\e}:
[a,b] \rightarrow B$, where $\epsilon$ varies over an open
interval containing $0$, such that the following conditions are
verified:
\begin{enumerate}
\item $\g_0=\g$;
\item for all $\e$, $\g_{\e}(a)=x$;
\item for any $\e \geq 0$ we have that
${\g}_{\e}$ is a concatenation of integral curves of vector fields
in ${\cal D}$;
\item the map $\e \mapsto \g_{\e}(b)$ is a smooth curve through $b$.
\end{enumerate}
The tangent vector to the curve $\e \mapsto \g_{\e}(b)$ at $\e =
0$ is called the {\em tangent vector to the variation $\g_{\e}$}
(note that $\g_0(b)=\g(b)= \Phi_T(x) =y$). Rather than considering
all possible variations satisfying the above conditions, we will
mainly deal with a specific class of variations, to be determined
below, called {\em single variations}. It will be shown that the
tangent vectors at $y$ to these single variations generate a
convex cone in $D_y$ (where we recall that $D$ refers to the
smallest integrable distribution generated by ${\cal D}$) and,
moreover, we will prove that each vector belonging to this cone is
in fact a tangent vector to a variation. If we agree to call {\em
dimension of a cone} the dimension of the linear space generated
by all vectors belonging to the cone, then the main result of this
section can be summarized as follows: if the dimension of the cone
of tangent vectors at $y$ to single variations equals the
dimension of $D_y$, say $d$, then there exists a coordinate chart
$V$ on the leaf $L_y$, with $y \in V$ and coordinate functions
denoted by $(x^1,\ldots,x^d)$, such that for any point $z \in V$
for which $x^i(z)\ge 0$ for all $i=1,\ldots,d$, we have that $z
\in R_x$.

Consider again a concatenation of integral curves $\g:[a,b]
\rightarrow B$ associated with the composite flow $\Phi: V \subset
\R^{\ell} \times B \rightarrow B$ of an ordered set of $\ell$
vector fields $(X_{\ell}, \ldots, X_1)$ in $\cal D$, and with a
given value $T \in \R^{\ell}_+$ of the corresponding composite
flow parameter, such that $\g(a)=x$ and $\g(b) = \Phi_T(x)=y$. We
now proceed towards the construction of what will be called a
single variation of $\g$. Let $T = (t_{\ell},\ldots,t_1) \in
\R^\ell_+$ and put $a_0=a,a_{\ell}=b$ and $ a_i= a_{i-1}+ t_i$ for
$i=1, \ldots, \ell$. Choose an arbitrary point $\t \in
]a_0,a_\ell]$ and let $Y$ be any vector field on $B$ such that
$\g(\t)$ belongs to the domain of $Y$. To fix the ideas, let us
assume that $a_{i-1} < \t \leq a_i$.  The flow of $Y$ will be
denoted by $\{\psi_s\}$ and, as before, $\{\phi^i_s\}$ denotes the
flow of $X_i$. We can then consider the composite flow $\Phi^* :
V' \subset \R^{\ell + 2} \times B \to B$, associated with the
ordered set of $\ell + 2$ vector fields $(X_\ell,\ldots,
X_i,Y,X_i, \ldots,X_1)$. Next, define
\begin{equation}\label{compospar}
\begin{array}{rcl}
T^*: \R& \rightarrow& \R^{\ell+2}:\\
\e& \mapsto &T^*(\e) =
(t_\ell,\ldots,t_{i+1},a_i-\t,\e,\t-a_{i-1},t_{i-1},\ldots,t_1)\,.
\end{array}
\end{equation}
It is easily seen that there exists an open neighborhood
$\tilde{I} \subset \R$ of $0$, such that $x$ is contained in the
domain of the map $\Phi^*_{T^*(\e)}$ for all $\e \in \tilde{I}$.
For each $\epsilon \in \tilde{I}$, let $\g_{\epsilon}$ denote the
concatenation of integral curves through $x$ corresponding to
$\Phi^*$ and $T^*(\e)$. The following sketch visualizes the
situation for $\t \in ]a_1,a_2]$:

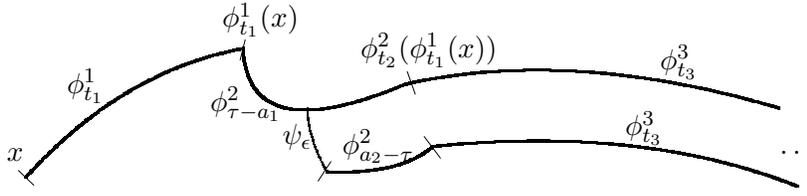
\begin{center}\setlength{\unitlength}{.8cm}
\begin{picture}(13.97,4.20) \thicklines
\qbezier(0.424,1.80)(2.01,3.56)(4.03,3.94)
\qbezier(4.03,3.94)(4.11,2.25)(6.74,3.35)
\qbezier(6.74,3.35)(9.51,3.96)(12.93,2.95) \thinlines
\drawline(0.318,1.91)(0.572,1.67)
\drawline(4.09,4.09)(4,3.75)
\put(3.7,4.3){{$\phi_{t_1}^{1}(x)$}}
\drawline(6.78,3.56)(6.91,3.18) 
\put(6,3.8){{$\phi_{t_2}^{2}(\phi_{t_1}^1(x))$}}
\put(1.12,3.2){$\phi_{t_1}^{1}$}
\put(3.5,2.9){{$\phi^{2}_{\t-a_1}$}}
\put(10.98,3.6){{$\phi_{t_3}^{3}$}} \put(0.127,2.08){\rmfamily
$x$} \put(12.97,2.2){\rmfamily \ldots} \thinlines
\qbezier[40](5.11,2.95)(5.11,2.46)(5.43,1.89) \thicklines
\qbezier(5.43,1.89)(6.65,1.84)(7.18,2.31)
\qbezier(7.18,2.31)(10.66,2.67)(13.24,1.65) \thinlines
\drawline(5.51,2.06)(5.30,1.74) \drawline(7.06,2.48)(7.33,2.12)
\put(5.7,2.2){$\phi^{2}_{a_2-\t}$}
\put(10.40,2.5){{$\phi_{t_3}^{3}$}} \put(4.7,2.35){\rmfamily
$\psi_\e$}
\end{picture}\end{center}

The tangent vector to the smooth curve $\e \mapsto
\g_{\e}(b)=\Phi^*_{T^*(\e)}(x)$ at $\e=0$ is then given by
\[
\left. \fpd{}{\e} \right|_{0} \Phi^*_{T^*(\e)} (x) = T
\Phi_{\t}^{a_\ell} (Y(\g(\t))) \in T_yB\,,
\]
where, in order to simplify the notations, we have introduced the
mapping $T\Phi_\t^{a_\ell}: T_{\g(\t)} B \to T_{y}B$, given by
\[
T\Phi_\t^{a_\ell}(v) = T \phi^\ell_{t_\ell} \circ T
\phi^{\ell-1}_{t_{\ell-1}} \circ \ldots \circ
T\phi^i_{a_i-\t}(v)\,, \quad \forall v \in T_{\g(\t)}B\,.
\]

Assume now that $Y \in {\cal D}$. Then one can see that the
1-parameter family of continuous piecewise curves $\g_{\e}$
satisfies the conditions proposed above for a variation of $\g$.

Next, suppose we take $Y = -X_i$ and $\t \in ]a_{i-1},a_i]$ for
some $i \in \{1, \ldots,\ell\}$, then for $\epsilon > 0$ (but
sufficiently small) and for any $t \in ]\t, \t+\epsilon]$, the
tangent vector $\dot{\g}_{\e}(t)$ to the concatenation of integral
curves through $x$, induced by $\Phi^*$ and $T^*(\e)$, in general
will not be contained in ${\cal D}$ since $-X_i$ does not have to
belong to ${\cal D}$. Consequently, if $-X_i \not\in {\cal D}$,
the $\g_{\e}$ resulting from the choice $Y=-X_i$ is, strictly
speaking, not a variation in the sense put forward above. However,
we can easily remedy the situation by constructing a ``reduced''
composite flow as follows. Putting $\widehat{T}(\e ) = (t_\ell,
\ldots, t_{i} - \e,\ldots, t_1) \in \R^{\ell}$, we see that for
$\epsilon$ sufficiently small, $\Phi_{\widehat{T}(\e)}$ is
well-defined in a neighborhood of $x$ and, moreover, since
$\phi^i_{a_i-\t} \circ \phi^i_{-\e} \circ \phi^i_{\t-a_{i-1}} =
\phi^i_{t_i -\e}$, it follows that $\Phi^*_{T^*(\e)} =
\Phi_{\widehat{T}(\e)}$. The concatenation of integral curves
determined by $\Phi$ and $\hat{T}(\e)$ does verify the conditions
for a variation of $\g$. The tangent vector at $\e =0$ to this
``reduced" variation equals
\[
\left. \fpd{}{\e} \right|_{0} \Phi^*_{T^*(\e)} (x) = \left.
\fpd{}{\e} \right|_{0} \Phi_{\widehat{T}(\e)}(x)=
-T\Phi^{a_{\ell}}_{\t}(X_i(c(\t)) ).
\]
We have thus shown that if $\t \in ]a_{i-1},a_i]$, a variation of
the given $\g$ is also determined by the ordered set
$(X_{\ell},\ldots,X_i,-X_i,X_i,\ldots,X_1)$.

To conclude, if we are given a continuous piecewise curve $\g:
[a,b] \rightarrow B$, with $\g(a)=x$, such that $\g$ consists of a
concatenation of integral curves determined by the composite flow
$\Phi$ and composite flow parameter $T = (t_{\ell},\ldots,t_1) \in
\R^{\ell}_{+}$ of an ordered set of vector fields
$(X_{\ell},\ldots,X_1)$ belonging to $\cal D$, we introduce the
following definition.
\begin{defn}
A {\rm single variation} of $\g$ is a $1$-parameter family of
continuous piecewise curves $\g_{\epsilon}:[a,b] \rightarrow B$,
passing through $x$, with $\g_0=\g$, and such that for each $\e$
the corresponding $\g_{\e}$ is the continuous piecewise curve
determined by the composite flow $\Phi^*$ and composite flow
parameter $T^*(\e)$ associated to an ordered set of vector fields
of the form $(X_{\ell},\ldots,X_i,Y,X_i,\ldots,X_1)$ for some $i
\in \{1,\ldots,\ell\}$, with $Y \in {\cal D} \cup \{-X_i\}$ and
where $T^*(\e)$ is given by (\ref{compospar}). (We will also
briefly refer to $\g_{\e}$ as `the single variation determined by
$\Phi^*$ and $T^*(\e)$'.)
\end{defn}

For later use we introduce the shorthand notation: ${\cal
D}_{-X}:= {\cal D} \cup \{-X_i\,|\,i=1,\ldots,{\ell}\}$. Whenever
we consider a single variation determined by an ordered set
$(X_{\ell},\ldots,X_i,Y,X_i,\ldots,X_1)$ for some $Y \in {\cal
D}_{-X}$, it will always be understood that $Y=-X_j$ can only
occur if $i=j$.

Given a single variation $\g_{\e}$ of $\g$, determined by a
composite flow $\Phi^*$ and composite flow parameter $T^*(\e)$,
one can always obtain a `new' variation by considering a suitable
reparameterization $\e(\e')$. More precisely, let $\e' \mapsto
\e(\e')$ denote a smooth map satisfying $\e(0) = 0$ and $\d =
\frac{d\e}{d\e'}(0)> 0$. Then it is not difficult to verify that
$\Phi^*$ and $T^*(\e(\e'))$ also determine a variation since
$\d>0$ implies that, in a neighborhood of $0$, $\sign(\e)=
\sign(\e')$. The tangent vector to the curve $\e' \mapsto
\Phi^*_{T^*(\e(\e'))}(m)$ at $\e'=0$ equals $\d \ T
\Phi_{\t}^{a_\ell} (Y(\g(\t)))$. From this one can easily derive
that any positive multiple of a tangent vector to a single
variation is again a tangent vector to a (not necessarily single)
variation. Note that if $\d Y \in {\cal D}_{-X}$, then $\d \ T
\Phi_{\t}^{a_\ell} (Y(\g(\t)))$ is again a tangent vector to a
single variation. In general, however, if $Y \in {\cal D}_X$, the
vector field $\d Y$ need not be contained in ${\cal D}_{-X}$. All
this naturally leads to the following definition.

\begin{defn}
Let $y \in R_x$ and fix a composite flow $\Phi$, corresponding to
an ordered set $(X_{\ell},\ldots,X_1)$ of vector fields in ${\cal
D}$, such that $\Phi_T(x)=y$ for some $T \in \R^\ell_+$. The {\rm
variational cone at $y$} associated to $\Phi$ and $T$, is the cone
$C_yR_x(\Phi,T)$ in $T_yB$ consisting of all finite linear
combinations, with positive coefficients, of tangent vectors to
single variations, i.e.
\begin{eqnarray*}
C_yR_x(\Phi,T) = \{&& \sum_{i=1}^s \d^i
T\Phi^{a_\ell}_{\t^i}(Y_i(\g(\t^i))) \ |\ Y_i \in {\cal D}_{-X},
\d^i \ge 0,\\ && \t^i \in ]a_0,a_\ell], s\in\N \}.
\end{eqnarray*}
\end{defn}
If no confusion can arise, we will often drop the explicit
reference to $\Phi$ and $T$ and simply denote the variational cone
by $C_yR_x$. It is easily seen that $C_yR_x$ is a convex set.
Indeed if $v,w \in C_yR_x$, then $(1-t)v + tw \in C_yR_x$, for any
$t \in [0,1]$. As a consequence of the next lemma it will be seen
that any element of $C_yR_x(\Phi,T)$ can be regarded as a tangent
vector to a variation of the continuous piecewise curve through
$x$ associated with $\Phi$ and $T$. First, we introduce an
alternative notation for composite flows which will sometimes be
more convenient, in particular when considering compositions of
composite flows.

Let $(Z_\ell,\ldots,Z_1)$ denote an ordered family of vector
fields on a manifold $B$, with composite flow $\Psi$. If
$\{\psi^i_s\}$ represents the flow of $Z_i$ for $i=1,\ldots,\ell$,
then it will turn out to be convenient to write $\psi^\ell \star
\ldots \star \psi^1$ for the composite flow $\Psi$, whereby it is
understood that $(\psi^\ell \star \ldots \star
\psi^1)_T:=\psi^\ell_{t_\ell} \circ \ldots \circ \psi^1_{t_1} =
\Psi_T$ for any admissible $T = (t_\ell,\ldots,t_1)$. Using this
notation, we are able to define the {\em composition} $\Psi_{(2)}
\star \Psi_{(1)}$ of two composite flows $\Psi_{(2)},\Psi_{(1)}$,
with $\Psi_{(i)} = \psi_{(i)}^{\ell_i} \star \ldots \star
\psi_{(i)}^1$ for $i=1,2$, as follows
\[
\Psi_{(2)} \star \Psi_{(1)} = \psi^{\ell_{2}}_{(2)} \star \ldots
\star \psi^1_{(2)} \star \psi^{\ell_1}_{(1)} \star \ldots \star
\psi^1_{(1)}.
\]

We now have the following result, the proof of which is quite
technical. As before, we start from a given continuous piecewise
curve $\g: [a_0,a_{\ell}] \rightarrow B$, with $\g(a_0)=x$,
associated to the composite flow of an ordered set of $\ell$
vector fields $(X_{\ell},\ldots,X_1)$ in ${\cal D}$, and a fixed
value $T$ of the composite flow parameter.
\begin{lem} \label{eindig}
Consider any finite number of (say, $s$) tangent vectors to single
variations of $\g$, namely $v_i =
T\Phi^{a_\ell}_{\t^i}(Y_i(\g(\t^i)))$, with $Y_i \in {\cal
D}_{-X}$ and $\t^i \in ]a_{0},a_\ell]$ for $i=1, \ldots, s$. Then,
there exists a composite flow $\Phi^*$ associated to $\ell +2s$
vector fields, and a smooth mapping $T^*: \R^s \to \R^{\ell +
2s}\,, (\e^1,\ldots,\e^s) \mapsto T^*(\e^1,\ldots,\e^s)$ such
that:
\begin{enumerate}
\item \label{eig1} $\Phi^*_{T^*(0)} = \Phi_T$; \item \label{eig2}
$x$ belongs to the domain of $\Phi^*_{T^*(\e^1, \ldots , \e^s)}$
for all $(\e^1, \ldots , \e^s)$ in some open neighborhood
$I^{(s)}$ of $(0,\ldots,0) \in \R^s$; \item \label{eig3} for each
fixed $(\e^1,\ldots,\e^s) \in I^{(s)}$, with $\e^i> 0$ for all
$i$, the tangent vector to the concatenation of integral curves
through $x$ determined by $\Phi^*$ and $T^*(\e^1,\dots ,\e^s)$ is
everywhere contained in ${\cal D}$ (possibly after a `reduction'
of $\Phi^*$ in the sense described above) such that, in
particular,$\Phi^*_{T^*(\e^1, \ldots , \e^s)}(x) \in R_x$; \item
\label{eig4} the tangent vector at $\e = 0$ to the curve $\e
\mapsto \Phi^*_{T^*(\e \d^1,\ldots,\e \d^s)}(x)$ equals $\d^i
v_i$, for all $\d^i \in \R$ (and where the curve is defined on a
sufficiently small interval such that $(\e\d^1,\ldots,\e\d^s) \in
I^{(s)}$).
\end{enumerate}
\end{lem}

\begin{pf}
Without loss of generality, we may assume that the instants $\t^i$
are ordered in such a way that $\t^1 \le \t^2 \ldots \le \t^{s}$.
Moreover, whenever some of the successive $\t^i$ coincide, the
ordering should be such that from the corresponding vector fields
$Y_i$, those that do not belong to ${\cal D}$ always precede those
that do belong to $\cal D$. More precisely, assume $\t^i = \ldots
=\t^j$ with $1 \leq i < j \leq s$, and let $\t^i \in
]a_{r-1},a_r]$ for some $r \in \{1,\ldots,\ell\}$. Then we require
that if $Y_k = -X_r$ for some $k \in \{i, \ldots, j\}$, and $-X_r
\not\in {\cal D}$, we have $k < k'$ for all those $k'\in
\{i,\ldots,j\}$ for which $Y_{k'} \in {\cal D}$. Such an
arrangement can always be achieved by simply taking a suitable
permutation of the ordered set $(Y_i, \ldots Y_j)$, if necessary.
Henceforth, we will always assume, for simplicity, that the
$Y_i$'s already appear in the correct ordering.

For $j=1,\ldots,\ell$, let $s_j$ denote the maximum of the set
$\{i\ |\ \t^i \in ]a_{j-1},a_j]\}$ and put $s_j = s_{j-1}$ if $\{i
\ | \ \t^i \in ]a_{j-1},a_j]\}= \emptyset$ and $s_0 = 0$. The
number of $\t^i$'s belonging to the $j$-th subinterval is then
given by $n_j = s_j - s_{j-1}$. Let $\{\psi^i_s\}$ denote the flow
of $Y_i$ (and, as before, $\{\phi^j_s\}$ refers to the flow of
$X_j$). Using the `star' notation introduced above, we now
consider for each $j \in \{1,\ldots,\ell\}$, the composite flow
$\Phi^*_j : \R^{1+2n_j} \times B \to B$ defined by
\[
\Phi^*_j = \left\{\begin{array}{lll}\phi^j\star \psi^{s_j}\star
\phi^j\star \psi^{s_j -1}\star \ldots\star \phi^{j}\star
\psi^{s_{j-1}+1}\star \phi^j &\mbox{ if }& n_j > 0,\\\phi^j
&\mbox{ if } & n_j =0,
\end{array}\right.
\]
and a mapping $T^*_j: \R^{n_j} \mapsto \R^{1+2n_j}$ (where it is
understood that if $n_j=0$, then $T^*_{j} \in \R$):
\[
T^*_j(\e^{s_{j-1}+1}, \ldots, \e^{s_j}) =
\left\{\begin{array}{l}(a_j-\t^{s_j}, \e^{s_j}, \t^{s_j}
-\t^{s_j-1}, \e^{s_j-1}, \ldots, \\ \ \ \t^{s_{j-1}+2} -
\t^{s_{j-1}+1} , \e^{s_{j-1}+1}, \t^{s_{j-1}+1} - a_{j-1})\\\mbox{
if } n_j >0,\\ (a_j-a_{j-1}) \mbox{ if } n_j =0.\end{array}\right.
\]
Next, by $\Phi^*$ we denote the `composition' of all the composite
flows $\Phi^*_j$, i.e. $\Phi^* = \Phi^*_{\ell} \star \ldots \star
\Phi^*_1$. Then, $\Phi^*$ itself is a composite flow which can be
evaluated at points of $\R^{\ell +2s}\times B$. If we define the
mapping $T^*: \R^s \rightarrow \R^{\ell + 2s}$ by
\[T^*(\e^1, \ldots,\e^s) = (T^*_\ell(\e^{s_{\ell-1}+1}, \ldots,
\e^{s_\ell}), \ldots, T^*_1(\e^{1}, \ldots , \e^{s_1})),
\]
then it is easily seen that $(T^*(0,\ldots ,0),x) \in
\mbox{Dom}(\Phi^*)$ and the equation $y=\Phi^*_{T^*(0,\ldots
,0)}(x)$ holds. This implies, in particular, that there exists an
open neighborhood $I^{(s)}$ of $(0,\ldots,0) \in \R^s$ for which
the map
$(\e^1,\ldots,\e^s)\mapsto\Phi^*_{T^*(\e^1,\ldots,\e^s)}(x)$ is
well defined and, hence, (\ref{eig2}) holds. Note that
$\Phi^*_{T^*(\e^1,\ldots,\e^s)}(x)$ can still be written as:
\[\Phi^*_{T^*(\e^1,\ldots, \e^s)}(x) =
(\Phi^*_\ell)_{T^*_\ell(\e^{s_{\ell-1}+1}, \ldots, \e^{s_{\ell}})}
\circ \ldots \circ (\Phi^*_1)_{T^*_1(\e^{1}, \ldots ,
\e^{s_1})}(x).\]

For $s=1$ the definitions of $\Phi^*$ and $T^*$ coincide with
those encountered in the construction of a single variation. For
any $(\d^1,\ldots,\d^s) \in \R^s$ and $\e$ varying over a
sufficiently small interval centered at $0$, such that the image
of the map $\e \mapsto (\e\d^1,\ldots, \e\d^s)$ is contained in
$I^{(s)}$, a straightforward, but rather tedious computation shows
that the tangent vector to the curve $\e \mapsto \Phi^*_{T^*(\e
\d^1, \ldots ,\e \d^s)}(x)$, at $\e =0$, equals $\d^i v_i$,
proving (\ref{eig4}). It is also easily seen that when putting
$\e^i=0$ for all $i$, we obtain $\Phi^*_{T^*(0)}= \Phi_T$, proving
(\ref{eig1}).

The proof of (\ref{eig3}) we will be provided for a particular,
simplified case from which the idea for the general proof can then
be easily deduced. Recall that we have chosen the ordering of the
$\t^i$ in such a way that, whenever we have a sequence
$\t^i,\ldots,\t^j, (i < j)$ with $\t^i=\t^{i+1}=\cdots =\t^j$,
those vector fields $Y_k$ which belong to the set $\{-X_1, \ldots
,-X_\ell\}$ and which are not contained in ${\cal D}$, always
appear before all the $Y_{k'} \in {\cal D}$ in the sequence
$Y_i,\ldots,Y_j$. Consider now the particular case where $a_0 <
\t^1=\t^2=\t^3 < a_1 < \t^4$, $Y_1= -X_1 (\not\in {\cal D})$ and
$Y_2,Y_3 \in {\cal D}$. Then,
\begin{eqnarray*}
(\Phi^*_1)_{T^*_1(\e^1,\e^2,\e^3)} &=& \phi^1_{a_1-\t^1} \circ
\psi^3_{\e^3} \circ \psi^2_{\e^2} \circ \phi^1_{-\e^1} \circ
\phi^1_{\t^1 - a_0}\\
&=&\phi^1_{a_1-\t^1} \circ \psi^3_{\e^3} \circ \psi^2_{\e^2} \circ
\phi^1_{\t^1-a_0-\e^1}\,. \end{eqnarray*} Therefore, we can define
a new composite flow, associated with vector fields in ${\cal D}$,
by putting $\widehat{\Phi}_1 = \phi^1\star \psi^3\star\psi^2\star
\phi^1$, and a new composite flow parameter
$\widehat{T}_1(\e^1,\e^2,\e^3) =
(a_{1}-\t^3,\e^3,\e^2,\t^1-a_0-\e^1)$. Then
$(\Phi_1)^*_{T_1^*(\e^1,\e^2,\e^3)} =
(\widehat{\Phi}_1)_{\widehat{T}_1(\e^1,\e^2,\e^3)}$ and, for
$\e^1$ sufficiently small, the components of
$\widehat{T}_1(\e^1,\e^2,\e^3)$ are positive, from which
(\ref{eig3}) readily follows for the `reduced' composite flow
$\widehat{\Phi}_1$ and the reduced composite flow parameter
$\widehat{T}_1$. A similar reasoning can be applied to the general
case, which completes the proof of the lemma.\qed
\end{pf}
The previous lemma implies, among others, that any $v$ in the cone
$C_yR_x(\Phi,T)$ can be regarded as a tangent vector to a
variation of the continuous piecewise curve $\g$ through $x$,
determined by $\Phi$ and $T$. Indeed, by definition of the cone
$C_yR_x(\Phi,T)$ we can always write $v$ (in a non-unique way) as
$v = \sum^s_{i=1}\d^i v_i$ for a finite number of tangent vectors
to single variations $v_i = T\Phi^{a_\ell}_{\t^i} Y_i(\g(\t^i))$,
with $\d^i > 0$. We can then associate to these $v_i$ a composite
flow $\Phi^*$, and a composite flow parameter
$T^*(\e^1,\ldots,\e^s)$, as in the above lemma. Then $\Phi^*$ and
$\e \mapsto T^*(\e \d^1,\ldots,\e \d^s)$ determine a one-parameter
family of continuous piecewise curves satisfying the conditions
for a variation of $\g$. Moreover, from the above lemma it follows
that the tangent vector to the curve $\e \mapsto \Phi^*_{T^*(\e
\d^1,\ldots,\e \d^s)}(x)$ at $\e = 0$ precisely equals $v$, which
we wanted to demonstrate.

Note that $C_yR_x(=C_yR_x(\Phi,T))$ is entirely contained in $D_y$
(with $D$, as before, the smallest generalized integrable
distribution generated by $\cal D$). If the dimension of the cone
$C_yR_x$ equals $d = \mbox{dim}\,D_y$, then this is equivalent to
saying that the the interior of the convex cone $C_yR_x$, with
respect to the standard vector space topology on $D_y$, is not
empty. Indeed, if we have $d$ independent vectors $v^1, \ldots
,v^d \in C_yR_x$, then the interior of the simplex in $D_y$,
determined by the ordered set $(0,v^1, \ldots,v^d)$, is contained
in $C_yR_x$. The converse is an immediate consequence of the fact
that any (nonempty) open ball in a vector space spans the full
space.

Before stating the main result of this section, we recall that
$L_y$ denotes the leaf of $D$ passing through $y$ (and, of course,
$L_y=L_x$). From the theory of integrable distributions, we know
that $L_y$ is an immersed submanifold of $B$ whose dimension
equals the rank of $D$ at $y$.
\begin{thm} \label{approx1}
Assume that the dimension of the cone $C_yR_x$ equals the
dimension $d$ of $D_y$. Then there exists a coordinate chart $V$
on the leaf $L_y$, with $y \in V$ and coordinate functions denoted
by $(x^1,\ldots,x^d)$, such that for any point $z \in V$ for which
$x^i(z)\ge 0$ for all $i=1,\ldots,d$, we have that $z \in R_x$.
\end{thm}
\begin{pf}
By assumption, the linear space spanned by all elements of
$C_yR_x$ equals $D_y$. We can therefore select a basis $\{v_1,
\ldots, v_d\}$ of the linear space $D_y$, with $v_i \in C_yR_x$
for all $i$. By definition of $C_yR_x$, each $v_i$ can then be
written as a finite linear combination of tangent vectors to
single variations, i.e. \begin{equation}\label{vi's} v_i =
\sum_{j=1}^{s_i}\d_{(i)}^j v_j^{(i)}, \;i=1,\ldots,d,
\end{equation}
for some $\d^j_{(i)}\in \R_+$, and where each $v_j^{(i)}$ is of
the form \[v_j^{(i)}=T\Phi^{a_\ell}_{\t_{(i)}^j}
Y^{(i)}_j(\g(\t_{(i)}^j))\] for some $Y_{j}^{(i)} \in {\cal
D}_{-X}$, $\t^j_{(i)} \in ]a_0,a_\ell]$. Although these
decompositions are not uniquely determined, for the remainder of
the proof we assume that for each of the given basis vectors $v_i$
one particular decomposition has been singled out, i.e. we make a
fixed choice for the $v_j^{(i)}$ and for the positive real numbers
$\d^j_{(i)}$ appearing in (\ref{vi's}). In total we thus have
$s=s_1+\ldots + s_d$ tangent vectors to single variations
$v_j^{(i)}$ which, however, need not all be different and/or
linearly independent. For convenience, we introduce the following
ordering:
$(v^{(1)}_1,\ldots,v^{(1)}_{s_1},v^{(2)}_1,\ldots,v^{(d)}_{s_d})$
and we denote an arbitrary element of this ordered set by
$w_{\a}$, with $\a = 1, \ldots s$ and such that $w_{\a} =
v_{\a}^{(1)}$ for $\a = 1, \ldots s_1$, $w_{\a}= v_{\a -
s_1}^{(2)}$ for $\a = s_1+1, \ldots, s_1 + s_2$, {\em etc.\ ...\
}. According to Lemma \ref{eindig} we can associate to the $s$
tangent vectors to single variations, $w_\a$, a composite flow
$\Phi^*$ and a map $T^*: \R^s \to \R^{\ell +2s}$ such that
\begin{enumerate}
\item $\Phi^*_{T^*(0)} = \Phi_T$, \item \label{eig6}
$\Phi^*_{T^*(\e^1, \ldots , \e^s)}(x) \in R_x$ if all $\e^i\ge 0$,
\item \label{eig5} for any fixed $(\d^1,\ldots,\d^s) \in \R^s$,
the tangent vector to the curve $\e \mapsto \Phi^*_{T^*(\e
\d^1,\ldots,\e \d^s)}(x)$ at $\e = 0$ equals $\d^\a w_\a$.
\end{enumerate}
With the convention that $s_0:=0$, we have for any $v \in D_y$
that
\[
v= l^i v_i = \sum_{i=1}^d\sum_{j=1}^{s_i} l^i\d^j_{(i)}
w_{s_0+\ldots +s_{i-1}+j} \in D_y\,.
\]
Putting \[ (\d^1,\ldots,\d^s) :=
(l^1\d^1_{(1)},\ldots,l^{1}\d^{s_1}_{(1)},l^2\d^1_{(2)},\ldots,l^d\d^{s_d}_{(d)})\,,
\]
we can still write $v$ as
\[ v = \sum_{\a=1}^{s}\d^{\a}w_{\a}\,.\]
Since the $\d^j_{(i)}$ in (\ref{vi's}) have been fixed, it follows
that all the coefficients $\d^{\a}$, appearing in this
decomposition of $v$, are determined unambiguously. Therefore, the
following mapping is well-defined:
\[ \widetilde{T}: D_y \to \R^{\ell+2s},\; v \mapsto \widetilde{T}(v)=
T^*(\d^1,\ldots,\d^s)\,,
\]
and, clearly, $\widetilde{T}$ is smooth.

From the properties of $\Phi^*$ and $T^*$, one can further deduce
that, on a sufficiently small open neighborhood $W$ of the origin
in the linear space $D_y$, the mapping given by
\[f: W (\subset D_y) \to B,\,v \mapsto
\Phi^*_{\widetilde{T}(v)}(x) \] is well-defined and smooth.
Moreover, by definition of $\Phi^*$, we have that $f(0) = y$ and
$\im f \subset L_y$. Let $j: L_y \hookrightarrow B$ denote the
natural inclusion and let us write $\widetilde{f}$ for $f$,
regarded as a mapping from $W$ into $L_y$, such that the following
relation holds: $j \circ \widetilde{f} = f$. Since $j$ is an
immersion and $f$ is smooth, it follows that $\widetilde{f}: W
(\subset D_y) \to L_y$ is smooth. In view of the natural
identification $T_0D_y \cong D_y$, it is easily proven, using
property (\ref{eig5}) of $\Phi^*$ and $T^*$, that the tangent map
of $f$ at $0$ satisfies, for any $v = \d^{\a}w_{\a} \in D_y$,
\[
T_0f(v) = \left.\frac{d}{d\e}\right|_0 f(\e v) =
\left.\frac{d}{d\e}\right|_0 \Phi^*_{T^*(\e\d^1,\ldots,\e\d^s)}(x)
= \d^\a w_\a =v\,.
\]
This, in turn, implies that $T_0\widetilde{f}: D_y \to
T_y(L_y)\equiv D_y$ is the identity map and, hence,
$\widetilde{f}$ induces a diffeomorphism from an open neighborhood
$\widetilde{W} \subset W$ of $0 \in D_y$ onto a an open
neighborhood $V$ of $y$ in $L_y$. Hence, to each point $z \in V$
there corresponds a unique $v \in \widetilde{W}$, with
$\widetilde{f}(v)=z$ and, with respect to the basis $\{v^i:\;
i=1,\ldots,d\}$ of $D_y$ chosen above, we can write $v = l^iv_i$.
The open set $V$ then becomes the domain of a local coordinate
chart on $L_y$, with coordinate functions $x^i\; (i=1,\ldots,d)$
defined by putting $x^i(z) = l^i$. Finally, from property
(\ref{eig6}) of $\Phi^*$ and $T^*$ it follows that for those
vectors $v=l^iv_i \in \widetilde{W}$ for which all $l^i \ge 0$, we
have $z=f(v) \in R_x$ since, in this case, all the coefficients
$\d^a$ appearing in the decomposition $v = \d^{\a}w_{\a}$ are also
non-negative. This completes the proof of the theorem. \qed
\end{pf}

Observe that the coordinate vector fields on $L_y$ corresponding
to the special chart constructed in the previous theorem are such
that (using the notations from the proof of the theorem)
$\displaystyle{\left.\frac{\partial}{\partial x^i}\right|_{y}=
T_0\widetilde{f}(v^i)=v^i}$. This observation will be of use in
proving the following result, which is a straightforward
consequence of Theorem \ref{approx1}.
\begin{cor} \label{approx2}
Assume that $C_yR_x$ has a non empty interior with respect to the
topology of $D_y$ (denoted by $\mbox{\rm int}(C_yR_x)$). Then, for
any curve $\theta: [0,1] \to (L_x =) L_y$ with $\theta(0)=y$ and
$0 \neq \dot \theta(0) \in \mbox{\rm int}(C_yR_x)$ there exists an
$\e
>0$ such that $\theta(t') \in R_x$ for $0\le t' \le \e$.
\end{cor}
\begin{pf}
As pointed out before, the fact that $C_yR_x$ has nonempty
interior implies that the `dimension' of the cone equals that of
$D_y$ and so the previous theorem applies. One can always fix a
basis $v_i$ in $D_y$, with $v_i \in C_yR_x$, such that the $\dot
\theta(0)$ is contained in the interior of the simplex spanned by
$(0,v^1,\ldots,v^d)$. In particular, this means that $\dot
\theta(0) = k^i v_i$ with all $k^i \in ]0,1[$. Consider the
coordinate chart $(x^1,\ldots,x^d)$ on $L_y$, in a neighborhood of
$y$, associated with the basis $v^1,\ldots,v^d$ as constructed in
Theorem \ref{approx1}. Note, in passing, that $x^i(y)=0$ for all
$i$. Now, since $\left.\fpd{}{x^i}\right|_y = v_i$ for
$i=1,\ldots,d$, and putting $\theta^i = x^i \circ \theta$, we find
that
\[
\left.\frac{d}{dt'}\right|_0 \theta^i(t')  = k^i, \mbox{ for }
i=1,\ldots,d.
\]
This implies that for all $i=1,\ldots,d$, $\dot{\theta}^i(0)>0$
and hence, since $\theta^i(0) = 0$, $\theta^i(t')> 0$ for $0\le t'
\le \e$ and $\e$ sufficiently small, i.e.\ $x^i(\theta(t')) >0$
for $i=1,\ldots,d$. According to Theorem \ref{approx1} this
implies that $\theta(t')\in R_x$ for all $0 \le t' \le \e$.\qed
\end{pf}

To close this section, we return to the framework of a geometric
control structure.

\subsection*{The vertical variational cone in a geometric
control structure}

Let $(\t,\nu,\ro)$ denote an arbitrary geometric control
structure. It is easily seen that the previous definitions and
results can be applied, in particular, to the everywhere defined
family of vector fields ${\cal D} = \{\T \circ \ro \circ \s \ | \
\s \in \G(\nu)\}$ on $M$. Consider a pair $(m,n) \in M \times M$
such that $m\stackrel{(\Phi,T)}{\longrightarrow} n$ and let
$C_nR_m(\Phi,T)$ denote the associated cone of variations. Since
$M$ is fibred over the real line, the kernel of the tangent map $T
\t$ defines a sub-bundle $V\t=\ker T\t$ of $TM$, called the
vertical bundle to $\t$. We will now define a `sub-cone' of
$C_nR_m$ which is vertical in the sense that it is contained in
$V_n\t$ and which satisfies $V_nR_m \subset C_nR_m$.

\begin{defn}\label{verticalcone}
The {\rm vertical variational cone at $n$}, associated to $\Phi$
and $T$, is given by: \begin{eqnarray*} V_nR_m(\Phi,T) = &&
\{\sum_{i=1}^s
\d^iT\Phi^{a_\ell}_{\t^i}(Y^i(c(\t^i)) - \dot c(\t^i))\\
 &&
 \hspace{1cm} | \ \d^i\ge0, \t^i \in ]a_{0},a_\ell], Y^i \in {\cal D},
i=1,\ldots,s \} \end{eqnarray*}
\end{defn}
As for the variational cone, we shall also sometimes simply write
$V_nR_m$ if there can be no confusion regarding the related $\Phi$
and $T$.
\section{The cost coordinate and optimality} \label{costoptimal}
In this section we give a straightforward application of Corollary
\ref{approx2} leading to necessary conditions to be satisfied by
an optimal control. We first specify how the notion of optimality
of a control can be formulated within the present geometric
framework.

Let $(\t,\nu,\ro)$ be an arbitrary geometric control structure
(with $\t: M \to \R$, $\nu:U \to M$, $\ro: U \to J^1\t$, as in
Definition 2.1) and let $L \in \cinfty{U}$ denote a function on
the control bundle $U$. If $u: I=[a,b] \to U$ is a control, then
the {\em cost} of $u$ with respect to $L$ is defined by
\[{\J}(u) =\int^{b}_{a} L(u(t))dt.\]
If we put $m=\nu(u(a))$ and $n=\nu(u(b))$, we have, with the
notations from Section \ref{sectiegeometricsetting}, that
$m\stackrel{u}{\to} n$ and, in particular, $n\in R_m$. We say that
the control $u$ is {\em optimal} if ${\J}(u) \le {\J}(u')$ for any
other control $u'$ such that $m \stackrel{u'}{\to} n$. For the
further discussion, it will be helpful to introduce the following
notation:
\[\J_u^{(t_1,t_2)} = \int^{t_2}_{t_1}L(u(t))dt,\] where $t_1,t_2
\in [a,b]$, with $t_1 \leq t_2$. Note that, in this notation,
${\J}(u) = {\J}_u^{(a,b)}$. The function $L$ is sometimes referred
to as the {\em cost function}.

\begin{defn}
A {\rm geometric optimal control structure} $(\t,\nu,\ro,L)$
consists of a geometric control structure $(\t,\nu,\ro)$ and a
cost function $L$.
\end{defn}

We will now show that to every geometric optimal control structure
$(\t,\ro,\nu,L)$ one can associate an {\em extended geometric
control structure}, $(\ovl \t,\ovl \nu,\ovl \ro)$ in which the
cost function is incorporated into the bundle map $\ovl\ro$. For
that purpose, we first introduce the product space $\ovl M: = M
\times \R$, the points of which will be denoted by $(m,J)$. For
reasons to become clear later on, $J$ will be called the {\em cost
coordinate}. The fibration $\t$ of $M$ over $\R$ induces the
fibration $\ovl \t: \ovl M \to \R\,, (m,J) \mapsto \ovl \t
(m,J)=\t(m)$. Next, for the extended control bundle we take $\ovl
U = U\times \R$, with projection onto $\ovl M$ given by $\ovl
\nu(u,J) = (\nu(u),J)$. Finally, we can define a bundle map $\ovl
\ro : \ovl U \to J^1\ovl \t$ as follows: $\ovl \ro(u,J) = (\ro(u),
J,L(u))$, where we have used the canonical identification between
$J^1\ovl \t$ and $J^1\t\times \R^2$ obtained as follows: given any
section $\ovl c(t) = (c(t),J(t))$ of $\ovl \t$, we map $j^1_t \ovl
c$ onto $(j^1_t c,J(t),\dot J(t))$. Note that $\ovl \t_{1,0}(\ovl
\ro(u,J)) = \ovl \nu(u,J)$ and, therefore, $(\ovl \t,\ovl \nu,\ovl
\ro)$ is indeed a well-defined geometric control structure.

Next, we shall prove that any control defined on a geometric
optimal control structure $(\t,\nu,\ro,L)$ induces a control on
the extended structure $(\ovl \t,\ovl \nu,\ovl \ro)$, and vice
versa. Let $u:I=[a,b]\to U$ be a control related to
$(\t,\nu,\ro,L)$, with $\nu(u(a))=m$ and $\nu(u(b)) =n$. We shall
construct a control $\ovl u$ in the associated structure $(\ovl
\t,\ovl \nu,\ovl \ro)$ such that for any $J_0 \in \R$ we have
$(m,J_0) \stackrel{{\ovl u}}{\to} (n,J_0 + {\J}_u^{(a,b)})$. More
precisely, define the map $\ovl u: I \to {\ovl U}$ by putting
\[\ovl u(t) = (u(t), J_0+{\J}_u^{(a,t)}).\]
It is easily seen that $\ovl u$ determines a piecewise section of
$\ovl \t\circ \ovl \nu$ whose projection onto $\ovl M$ is a
continuous piecewise section. Furthermore, the first-order jet of
the base section $\ovl\nu \circ \ovl u$ equals $j^1_t(\nu\circ
u,J_0+{\J}_u^{(a,t)}) = (j^1_t(\nu \circ u), J_0+{\J}_u^{(a,t)},
L(u(t)))$. Since $u$ is a control, we readily obtain the equality
$\ovl \ro \circ \ovl u = j^1 (\ovl \nu \circ \ovl u)$, which
implies that $\ovl u$ is indeed a control. On the other hand, the
projections of $\ovl u(a)$ and $\ovl u(b)$ onto $\ovl M$ are given
by $(m,J_0)$ and $(n,J_0+{\J}(u))$, respectively. It follows that
$(m,J_0) \stackrel{\ovl u}{\to} (n,J_0+{\J}(u))$ for the extended
geometric control problem (and for arbitrary $J_0 \in \R$).

Conversely, let $\ovl u: [a,b] \to \ovl U,\, t \mapsto \ovl u(t) =
(u(t),J(t))$ represent a control on the extended geometric control
structure $(\ovl \t,\ovl \nu,\ovl \ro)$. Then, if the base section
is written as $(\ovl \nu \circ \ovl u)(t) = \ovl c(t) =
(c(t),J(t))$ we can deduce from $\ovl \ro \circ \ovl u = j^1\ovl
c$ that $(\ro \circ u)(t) = j^1_t c$, i.e.\ $u: [a,b] \to U$ is a
control. Moreover, the cost coordinate satisfies $\dot J(t) =
L(u(t))$ and, hence,
\[
J(t) = J(a) + \int^t_{a}L(u(t)) dt=J(a) + {\J}_u^{(a,t)},
\]
In particular, we have $J(b) =J(a) + {\J}_u^{(a,b)}$.

Summarizing the preceding discussion, we have proven the following
result.

\begin{prop}
Let $(\t,\nu,\ro,L)$ denote a geometric optimal control structure.
Then for any $m,n \in M$ and $J_m,J_n \in \R$, we have that $m
\stackrel{u}{\to} n$ and ${\J}(u)=J_n-J_m$ for some control $u$
iff  $(m,J_m) \stackrel{(u,J)}{\to} (n,J_n)$ in the associated
extended geometric control structure, where $J:[a,b] \to \R$ is
given by $J(t) = J_m + {\J}_u^{(a,t)}$.
\end{prop}
Consider once more an arbitrary geometric optimal control
structure $(\t,\nu,\ro,L)$ and assume $m \stackrel{u}{\to} n$ for
some control $u$. According to the previous proposition we then
know that, for any $J_0 \in \R$, one can define an appropriate
function $J(t)$ such that $(m,J_0) \stackrel{(u,J)}{\to}
(n,J_0+{\J}(u))$. Let $\ovl c = \ovl \nu \circ (u,J)$ be the base
section of the control $(u,J)$. On $\ovl M$ we can then consider
the variational cone $C_{(n,J_0+J(u))}R_{(m,J_0)}$, resp.\ the
vertical variational cone $V_{(n,J_0+J(u))}R_{(m,J_0)}$,
associated to a composite flow $\ovl \Phi$ and composite flow
parameter $\ovl T$ determining the controlled section $\ovl c$,
with $\ovl \Phi_{\ovl T}(m,J_0) = (n,J_0+{\J}(u))$. The proof of
the following proposition relies on Corollary \ref{approx2}.

\begin{prop}
Let $(\t,\nu,\ro,L)$ denote a geometric optimal control structure
and assume $m \stackrel{u}{\to} n$ for a control $u$ which is
optimal. Then the interior of $C_{(n,J_0+{\J}(u))}R_{(m,J_0)}$
does not contain the tangent vector
$\left.-\fpd{}{J}\right|_{(n,J_0+{\J}(u))}$.
\end{prop}
\begin{pf} Assume that $(-{\partial}/{\partial J})_{(n,J_0+{\J}(u))}
\in \mbox{\rm int}(C_{(n,J_0+{\J}(u))}R_{(m,J_0)})$. Consider the
`vertical' curve $\theta(t) = (n,J_0+{\J}(u) -t)$ in $\ovl M$,
defined for $t\in[0,1]$, whose tangent vector at $t=0$ precisely
equals $(-{\partial}/{\partial J})_{(n,J_0+{\J}(u))}$. From
Corollary \ref{approx2} it then follows that there exists an $\e
>0$, sufficiently small, such that $\theta(t) \in R_{(m,J_0)}$ for
$t \in [0,\e]$. From this, one can deduce that there exists a
control $\ovl u'$ for which $(m,J_0) \stackrel{\ovl u'}{\to}(n,
J_0 + {\J}(u) - \e)$. In view of previous considerations, this
further implies that there exists a control $u'$ on
$(\t,\nu,\ro,L)$ such that $m\stackrel{u'}{\to} n$, with cost
${\J}(u') = {\J}(u)-\e$, ${\J}(u')<{\J}(u)$. Since $u$ was assumed
to be optimal, this clearly leads to a contradiction. \qed
\end{pf}
Before proceeding, we first recall some properties and terminology
regarding linear spaces and convex cones in a linear space. Let
${\cal V}$ be an arbitrary (finite dimensional) linear space. A
hyperplane in ${\cal V}$ (i.e. a linear subspace of co-dimension
one) can always be defined as the set of all vectors $v \in {\cal
V}$ satisfying $\langle \eta,v\rangle = 0$ for some (non-zero)
co-vector $\eta \in {\cal V}^*$. Such a hyperplane divides ${\cal
V}$ into two `half-spaces' which are given by the set of all $v$
such that $\langle \eta,v \rangle \le 0$, resp. $\langle
\eta,v\rangle \ge 0$, and which are called the `negative'
half-space and the `positive' half-space, respectively. If $C$ is
a convex cone in ${\cal V}$ which does not span the full space,
then there always exists a hyperplane such that $C$ is contained
in one of the corresponding half-spaces.

If we now return to the situation described in the previous
proposition, it follows from the above considerations that, under
the conditions of Proposition 4.3, there exists a hyperplane in
the tangent space $T_{(n,J_0 + {\J}(u))}\ovl M$ such that the
variational cone $C_{(n,J_0 + {\J}(u))}R_{(m,J_0)}$ is contained
in, say, the corresponding negative half-plane, whereas the vector
$(-{\partial}/{\partial J})_{(n,J_0+{\J}(u))}$ belongs to the
positive half-plane. From the fact that the vertical variational
cone $V_{(n,J_0+{\J}(u))}R_{(m,J_0)}$ is a subset of
$C_{(n,J_0+J(u))}R_{(m,J_0)}$, contained in the vertical subspace
$V_{(n,J_0+{\J}(u))}\ovl \t$, the following result is a
straightforward consequence of Proposition 4.3.
\begin{cor}\label{maximum1}
If $m\stackrel{u}{\to}n$ and if $u$ is optimal, then there exists
a hyperplane in $V_{(n,J_0 + {\J}(u))}\ovl \t$, determined by some
$\ovl\eta \in V_{(n,J_0+J(u))}^*\ovl \t$ (the dual space of the
vertical tangent space $V_{(n,J_0+{\J}(u))}\ovl \t$) such that
\begin{enumerate}
\item $\langle \ovl\eta , -\left.\fpd{}{J}\right|_{(n,J_0+{\J}(u))}
\rangle \ge 0$, and
\item $\langle \ovl\eta, v  \rangle\le 0$ for all
$v \in V_{(n,J_0+{\J}(u))}R_{(m,J_0)}$.
\end{enumerate}
\end{cor}

In order to relate the previous result to a more familiar
formulation of the necessary conditions for an optimal control, in
terms of solutions of differential equations, we will need a minor
generalization of the theory of connections over a bundle map as
developed, for instance, in \cite{algcon}.

\section{Lifts over bundle maps} \label{sectie5}
For the sake of completeness, we first briefly recall the setting
for defining {\em a lift over a bundle map}.

Consider a smooth manifold $B$ and a fibre bundle $\nu :N \to B$,
equipped with a bundle map $\an: N \to TB$ fibred over the
identity, as shown in the following commutative diagram.

\setlength{\unitlength}{.8cm}
\begin{picture}(1,4)(-3.5,2)
\thicklines \put(2,2.5){$B$} \put(2.7,2.7){\vector(1,0){2.5}}
\put(5.5,2.5){$B$} \put(5.7,4.3){\vector(0,-1){1.3}}
\put(5.5,4.6){$TB$} \put(2.7,4.8){\vector(1,0){2.5}}
\put(2,4.6){$N$} \put(2.2,4.3){\vector(0,-1){1.3}}
\put(3.8,2.2){$\een_B$} \put(3.8,5.2){$\an$} \put(6,3.5){$\t_B$}
\put(1.7,3.5){$\nu$}
\end{picture}

Note that, unlike the treatment in \cite{algcon}, we do not
require $N$ to be a vector bundle. Next, let $\pi: E \to B$ denote
an arbitrary fibre bundle over $B$ and consider the pull-back
bundle $\pi^*N$. We can then define the following notion of lift.

\begin{defn} A {\rm lift over} $\an$ is a bundle map $h:\pi^*N \to TE$
fibred over the identity on $E$ such that the following diagram
commutes:

\setlength{\unitlength}{.8cm}
\begin{picture}(13,4)(-3.5,2)
\thicklines \put(2,2.5){$N$} \put(2.7,2.7){\vector(1,0){2.5}}
\put(5.5,2.5){$TB$} \put(5.7,4.3){\vector(0,-1){1.3}}
\put(5.5,4.6){$TE$} \put(2.7,4.8){\vector(1,0){2.5}}
\put(1.7,4.6){${\pi}^{\ast}N$} \put(2.2,4.3){\vector(0,-1){1.3}}
\put(3.8,2.1){$\an$} \put(3.8,5.2){$h$} \put(6,3.5){$T{\pi}$}
\put(1.5,3.5){$\tilde{\pi}_2$}
\end{picture}
\end{defn}
A lift $h$ over $\an$ allows us to define the {\em $h$-lift} of a
section $s$ of $\nu$. More precisely, the $h$-lift of $s \in
\G(\nu)$ is a section of $\t_E$ defined by $s^h(e) =
h(e,s(\pi(e)))$, for all $e \in E$. Note that $s^h$ determines a
vector field on $E$.

A {\em $\an$-admissible curve} $c:I=[a,b]\to N$ is a smooth curve
such that the base curve $\nu \circ c = \tilde c$ in $B$ satisfies
$\dot{\tilde c}(t) =\an(c(t))$. If we assume that $\an(n) \neq 0$
for all $n \in N$, then any $\an$-admissible curve is a
concatenation of integral curves of vector fields belonging to the
family ${\cal D'} = \{\an \circ s \ | \ s \in \G(\nu)\}$. Indeed,
let $c: I \to N$ denote a $\an$-admissible curve, with base curve
$\tilde c$. Then $\dot {\tilde c}(t) \neq 0$ for all $t$, i.e.\
$\tilde c$ is an immersion. Following an argument of S. Helgason
(see \cite[p 28]{helgason}), one can prove that there exists a
finite subdivision $\{I_i\}$ of $I$ such that for the restriction
of $c$ to each of these subintervals $I_i$ there exists a local
section $s_i$ of $\nu$ verifying $s_i(\tilde c(t)) = c(t)$ for all
$t \in I_i$. It is easily seen that ${\tilde c}_{|I_i}$ is an
integral curve of $\an \circ s_i$.

\begin{rem}\textnormal{We can apply all this to a geometric control
structure $(\t,\nu,\ro)$, where we take $B=M, N=U, \an = \T \circ
\ro$. A control can then be equivalently characterized as a $(\T
\circ \ro)$-admissible curve $u: I \to U$, with the additional
constraint that it should be a section of $\t \circ \nu$, i.e.$(\t
\circ \nu)(t)=t$ for all $t$. We also recover here the property
that each $(\T \circ \ro)$-admissible curve is a concatenation of
integral curves of vector fields in ${\cal D}$.}
\end{rem}

Assume now that the bundle $E$ is a vector bundle and let $\Delta$
be the dilation vector field on $E$, with flow $\{\d_t\}$. A lift
$h$ over $\an$ is then said to be {\em linear} if $T\d_t \circ h
(n,e) = h(n,\d_t(e))$ for any $t$. Consider bundle adapted
coordinate charts on $N$ and $E$, denoted by $(x^i,n^\a)$ and
$(x^i,e^A)$, respectively. In coordinates, $h$ then reads
\[
h(x^i,n^\a,e^A) = \an^j(x^i,n^\a)\left.\fpd{}{x^j}\right|_e +
\G^A(x^i,n^\a,e^A) \left.\fpd{}{e^A}\right|_e,
\]
and $h$ is a linear lift iff $\G^A(x^i,n^\a,e^A) =
\G^A_B(x^i,n^\a)e^B$. The functions $\G^A_B$ are called the {\em
coefficients of $h$}. For the remainder of this section, we always
take $E$ to be a vector bundle (over $B$).

Given a linear lift $h$ and a $\an$-admissible curve $c: [a,b] \to
N$, with base curve $\tilde {c}$, take any $e \in E$ such that
$\pi(e) = \tilde c(a)$. We can then construct a curve $c^h$ in $E$
through $e$, called the {\em $h$-lift} of $c$, which is uniquely
determined by the differential equation $h(c^h(t),c(t)) = \dot
c^h(t)$, with initial condition $c^h(a)=e$ (see also
\cite{algcon}).

Next, we show that a linear lift $h$ always induces a derivative
operator $\del$, acting on sections of $\pi$. Let $\pi_2: V\pi
\cong E \times_B E \to E$ denote the projection onto the second
factor, then, in analogy with the case where $N$ is a vector
bundle and $\an$ a linear bundle map (see \cite{algcon}), we can
define a mapping $K: \an^*TE \to E$ according to: $K(n, w) =
\pi_2(w-h(\t_E(w),n))$. Given any $x \in B$, $n \in N_x (=
\nu^{-1}(x))$ and any local section $\psi \in \G(\pi)$, defined on
an open neighborhood of $x$, we put
\[ \del_n \psi:= K(n,T_x\psi(\an(n))).\]
Clearly, $\del_n\psi \in E_x (= \pi^{-1}(x))$. The map $\del_n$
thus defined, is a derivative operator on $\G(\pi)$ since, for
arbitrary $f \in \cinfty{B}$, $\psi_1, \psi_2 \in \G(\pi)$ (all at
least defined on a neighborhood of $x$) we find
that \begin{eqnarray*} &&\del_n f\psi = \an(n)(f) \psi(x) + f(x)\del_n \psi,\\
&&\del_n (\psi_1 + \psi_2) = \del_n \psi_1 + \del_n
\psi_2.\end{eqnarray*} An operator on $\G(\pi)$ satisfying these
properties is called a {\em $\an$-derivative}. Given any section
$s \in \G(\nu)$, we can define the operator $\del_s$ on $\G(\pi)$
by
\[
\del_s\psi(x):=\del_{s(x)}\psi,
\]
and, obviously, $\del_s\psi$ is again a section of $\pi$. It is
easily seen that there is a one-to-one correspondence between
$\an$-derivatives and linear lifts over $\an$. Using the above
coordinate expression for $h$, we obtain that the $\an$-derivative
determined by $h$ locally reads (for $n=(x^i,n^{\a}) \in N_x$)
\[(\del_n \psi)^A = \an^j(x^i,n^\a)\fpd{\psi^A}{x^j}(x^i) -
\G^A_B(x^i,n^\a)\psi^B(x^i).\] It also follows that $\del_s\psi=0$
for $s \in \G(\nu)$ and $\psi \in \G(\pi)$ iff
\[T_x\psi(\an(s(x))) = s^h (\psi(x))\] for all $x \in B$.

Similar to what we have in standard connection theory, a
derivative operator can be constructed which acts on sections of
$\pi$ defined along the base curve of a $\an$-admissible curve
$c:I=[a,b] \to N$. Indeed, consider a curve in $E$, $\tilde \psi :
I \to E$, such that $\pi \circ \tilde\psi = \nu \circ c (= \tilde
c)$, then the $\an$-derivative associated to the linear lift $h$
and acting on $\tilde \psi$ equals
\[
\del_c \tilde \psi(t) := K(c(t), \dot{\widetilde \psi}(t)).
\]
It is not difficult to prove that $\del_c \tilde \psi(t)=0$ for
all $t$ iff $\tilde \psi =c^h$. If $\del_c \tilde \psi \equiv 0$,
we say that $\tilde \psi$ is {\em $h$-transported along $c$} and
that $\tilde \psi(b)$ is the {\em $h$-transport} of $\tilde
\psi(a)$ along $c$. We conclude by pointing out that any
$\an$-admissible curve $c$ in $N$ determines a linear map
$c^{b}_a: E_{c(a)}\to E_{c(b)}$, called the {\em $h$-transport
operator along $c$}, defined by $c^b_a(e) = \tilde \psi(b)$, where
$\tilde \psi$ is the unique solution of the equation $\del_c\tilde
\psi (t)=0$ with $\tilde \psi(a)=e$.

\section{The control lift and control derivative}
\label{sectiecontrolafgeleide} Let $(\t,\nu,\ro)$ denote a
geometric control structure. Consider the first-order jet bundle
$J^1\nu$ of the bundle $\nu: U \to M$, with associated projections
$\nu_1: J^1\nu \to M$, $\nu_{1,0}: J^1\nu \to U$. Recall that for
any two local sections $\s$ and $\s'$ of $\nu$, defined on a
neighborhood of a point $m \in M$, we have that $j^1_m\s =
j^1_m\s' \in J^1\nu$ iff $\s(m) = \s'(m)$ and $T_m\s = T_m\s'$ (as
linear maps from $T_mM$ into $T_{\s(m)}U$). Bearing this in mind,
it is easily seen that the following mapping is well-defined:
\begin{equation}\label{lambda}
\an: J^1\nu \to TU,\, j^1_m\s \mapsto \an(j^1_m\s) =
T_m\s((\T\circ \ro)(\s(m))). \end{equation} Moreover, $\an$ is a
bundle map over the identity on $U$. In terms of appropriate
bundle coordinates $(t,x^i,u^a)$ on $U$ and
$(t,x^i,u^a,u^a_t,u^a_i)$ on $J^1\nu$, $\an$ reads
\[
\an(t,x^i,u^a,u^a_t,u^a_i) =
(t,x^i,u^a,1,\ro^j(t,x^i,u^a),u^b_t+\ro^j(t,x^i,u^a)u^b_j)).
\]
We now consider the fibred product bundle $U \times_M V\t$, with
projections $p_1:U\times_M V\t \to U,\, (u,v) \mapsto p_1(u,v) =
u$ and $p_2: U \times_M V\t \to V\t,\, (u,v) \mapsto v$, whereby
$\nu \circ p_1=\t_M \circ p_2$. Observing that $p_1: U \times_M
V\t \to U$ is a vector bundle over $U$, we can apply the theory
from the previous section to the case where $B = U$, $N = J^1\nu$,
$E = U \times_M V\t$ and $\an$ is given by (\ref{lambda}). It will
be seen that, within this setting, $\an$-admissible curves are
closely related to controls. For that purpose, we need the
following straightforward extension of the definition of
$\an$-admissible curve to the class of piecewise curves: a
piecewise curve $\psi$ in $J^1\nu$ is said to be $\an$-admissible
if it is induced by (i.e.\ consists of a concatenation of) a
finite number of smooth $\an$-admissible curves.

In the sequel, we always assume that a piecewise $\an$-admissible
curve $\psi$ in $J^1\nu$ has a continuous projection onto $M$ and
is parameterized such that $\t(\nu_1(\psi(t)))=t$, i.e. such that
$\psi$ is a section of $\t \circ \nu_1$. (Note that this is not a
restriction since, given any $\an$-admissible curve $\psi:[a,b]
\to J^1\nu$, with $t_a=\nu_1(\psi(a))$, we can consider a
reparametrization of $\psi$ according to $\psi':[t_a,t_a+b-a] \to
J^1\nu,\, t \mapsto \psi'(t)= \psi(t-t_a+a)$. Then, $\psi'$ is
still $\an$-admissible and, moreover, satisfies
$\t(\nu_1(\psi'(t)))=t$.)
\begin{lem} \label{lemma21}
The projection onto $U$ of any $\an$-admissible curve in $J^1\nu$
is a smooth control, and any control in $U$ can be obtained as the
projection of a piecewise $\an$-admissible curve.\end{lem}
\begin{pf}
We first prove that the projection $u = \nu_{1,0}\circ \psi$ of a
$\an$-admissible curve $\psi: [a,b] \to J^1\nu$ is a smooth
control. By definition, we have $\dot u(t) = \an(\psi(t))$. From
$T\nu \circ \an = \T \circ \ro \circ \nu_{1,0}$, it follows that
$\dot c(t) = (\T \circ \ro) (u(t))$, where $c = \nu_1(\psi(t)) =
\nu \circ \nu_{1,0}(\psi(t))$. This shows that the smooth curve
$u$ is $(\T \circ \ro)$-admissible, i.e. it is a smooth control.

On the other hand, assume that $u: [a_0,a_\ell] \to U$ is a
control, with base curve $c = \nu \circ u$. We then know that $c$
can be written as a concatenation of integral curves, induced by
the composite flow $\Phi$ of an ordered set $(\T \circ \ro \circ
\s_\ell, \ldots, \T \circ \ro \circ \s_1)$ for some $\s_i \in
\G(\nu)$, with composite flow parameter $T=(a_\ell -a_{\ell-1},
\ldots, a_1-a_0)$. Furthermore, $u(t) = \s_i(c(t))$ for any $t\in
]a_{i-1},a_{i}]$. Putting $\psi_i(t) = j^1\s_i(c(t))$ for all
$t\in [a_{i-1},a_{i}]$ and $i=1,\ldots,\ell$, we obtain that for
any $t \in ]a_{i-1},a_i]$ the equality
\[\an(\psi_i(t)) = T_{c(t)}\s_i(\dot c(t)) =
\left.\frac{d}{dt}\right|_t(\s_i(c(t)))= \dot u(t),\]holds.
Therefore, according to the definition above, the piecewise curve
$\psi: [a_0,a_\ell] \to J^1\nu$, induced by the smooth curves
$\psi_i:[a_{i-1},a_i]\to J^1\nu$, is a piecewise $\an$-admissible
curve, which completes the proof of the lemma.\qed
\end{pf}
In the following we shall frequently make use of the natural
identification $T(U \times_M V\t) \cong TU \times_{TM} T(V\t)$,
without mentioning it explicitly. We further denote by $\mbox{\fr
s}:TTM \to TTM$ the canonical involution on $TTM$. The latter is
characterized by the relations $T\t_M \circ \mbox{\fr s} =
\t_{TM}$ and $\t_{TM}\circ \mbox{\fr s} = T\t_M$.

\begin{rem}\textnormal{Recall that, given an arbitrary manifold $B$
with local coordinates $(x^i)$, and denoting the natural bundle
coordinates on $TB$ and $TTB$ by $(x^i,v^i)$ and
$(x^i,v^i,\dot{x}^i,\dot{v}^i)$, respectively, then the canonical
involution $\mbox{\fr s}$ on $TTB$ reads $\mbox{\fr s}\,
(x^i,v^i,\dot{x}^i,\dot{v}^i)=(x^i,\dot{x}^i,v^i,\dot{v}^i)$.}
\end{rem}

For a geometric control structure $(\t,\nu,\ro)$, with bundle map
$\an$ given by (\ref{lambda}), we have the following property.
\begin{prop}
The map $h^c : \nu_{1,0}^*(U\times_M V\t ) \to T(U \times_M V\t)$,
defined by \[h^c(j^1_m\s,(\s(m),v)) = (\an(j^1_m\s), \mbox{\fr
s}\,(T (\T \circ \ro)(T_m\s(v)))),\]for any $m \in M$, $\s \in
\G(\nu)$ and $v \in V_m\t (\subset T_mM)$, is a linear lift over
$\an$.
\end{prop}
\begin{pf}
We first verify that $h^c$ indeed takes values in $T(U \times_M
V\t)$. For that purpose, consider bundle adapted coordinates
$(t,x^i,u^a)$ and $(t,x^i,v^j)$ on $U$ and $V\t$, respectively.
Take $m=(t,x^i) \in M$, $v=(t,x^i,v^j) \in V_m\t$ and $j^1_m\s =
(t,x^i,u^a,u^a_t,u^a_i) \in J^1\nu$, then: \begin{eqnarray*}
&&\mbox{\fr s}\,(T (\T \circ \ro)(T_m\s(v))) =
\left.\fpd{}{t}\right|_{v} +
\ro^i(t,x^j,u^a)\left.\fpd{}{x^i}\right|_v \\ && +
\left(v^i\fpd{\ro^k}{x^i}(t,x^j,u^a) +
v^iu^b_i\fpd{\ro^k}{u^b}(t,x^j,u^a)\right)\left.\fpd{}{v^k}\right|_v.
\end{eqnarray*} From this expression one can read that \[\mbox{\fr s}\,(T (\T
\circ \ro)(T_m\s(X))) \in T_v(V\t)(\subset T_vTM).\] Next, using
the properties of the canonical involution operator $\mbox{\fr
s}$, and taking into account (\ref{lambda}), it is easily seen
that \[T\nu(\an(j^1_m\s))= T\t_M(\mbox{\fr s}\,(T (\T \circ
\ro)(T_m\s(v)))) =(\T \circ \ro)(\s(m)) \in T_mM\] which proves
indeed that $\im h^c \subset T(U \times_M V\t)$.

From its definition it readily follows that $h^c$ is a bundle map
fibred over the identity on $U\times_M V\t$, and we have that \[T
p_1 (h^c(j^1_m\s , (\s(m),X))) = \an(j^1_m\s).\] This already
guaranties that $h^c$ is a lift over $\an$ in the sense of
Definition 5.1. From the above coordinate expression we can also
deduce that the ${\partial}/{\partial v^k}$-components $\G^k$ of
$h^c$ are linear in the fibre coordinates $v^i$ of the vector
bundle $p_1$. More precisely, we have $\G^k(t,x^j,v^j,u^a,u^a_j) =
\G^k_i(t,x^j,u^a,u^a_j)v^i$ with
\[\G^k_i(t,x^i,u^a,u^a_i) = \fpd{\ro^k}{x^i}(t,x^j,u^a) +
u^b_i\fpd{\ro^k}{u^b}(t,x^j,u^a).\] This shows, in particular,
that $h^c$ is a linear lift. \qed
\end{pf}
Note that the `coefficients' $\G^k_i$ of $h^c$ do not depend on
the coordinates $u^a_t$ of $J^1\nu$. In a remark at the end of
this section we will return to this point in more detail.

Let us denote the $\an$-derivative corresponding to $h^c$ by $D$
and let $\vvectorfields{\nu}$ denote the set of $\t$-vertical
vector fields along $\nu$, i.e. \[\vvectorfields{\nu}= \{Z: U \to
V\t\,|\,\t_M(Z(u)) = \nu(u),\; \mbox{for all}\;u\in U\}.\] Note
that, in view of the relation $\nu \circ p_1 = \t_M \circ p_2$, we
have $\vvectorfields{\nu} \cong \G(p_1)$.
\begin{prop}
Given any $j^1_m\s \in J^1\nu$ and $Z \in \vvectorfields{\nu}$,
then $D_{j^1_m\s} Z$ is contained in $V_m\t$ and \[D_{j^1_m\s} Z =
[\T \circ \ro \circ \s , Z \circ \s](m),\] (where the square
brackets on the right-hand side denote the ordinary Lie bracket of
vector fields on $M$).
\end{prop}
\begin{pf}
Recalling the coordinate expression of a $\an$-derivative (cf.\
Section 5), we obtain, with a slight abuse of notation,
\[
(D_{j^1_m\s} Z)^i = (\fpd{Z^i}{t} + \ro^j \fpd{Z^i}{x^j} + (u^a_t+
\ro^j u^a_j) \fpd{Z^i}{u^a} - \G^i_j Z^j)_m.
\]
The result then easily follows upon substituting $u^a_j =
\fpd{\s^a}{x^j}$ and $u^a_t = \fpd{\s^a}{t}$ in the right-hand
side, and comparing this with the coordinate expression of the Lie
bracket $[\T \circ \ro \circ \s , Z \circ \s](m)$. \qed
\end{pf}
We shall now derive an explicit expression for the $h^c$-transport
operator $\psi^b_a$ determined by a $\an$-admissible curve $\psi:
I=[a,b] \to J^1\nu$. We first consider the case where $\psi$ takes
the special form $\psi(t)=j^1\s(c(t))$ for some curve $c:[a,b] \to
M$ and a section $\s \in \G(\nu)$. Note that such a $\psi$ is
$\an$-admissible iff $u(t) := \s(c(t))$ is a smooth control, which
still implies that in terms of the flow $\{\phi_s\}$ of the vector
field $\T \circ \ro \circ \s$, we have $c(t) = \phi_{t-a}(c(a))$.
\begin{lem}
Let $\psi: [a,b] \to J^1\nu,\, t \mapsto \psi(t)=j^1\s(c(t))$ be a
$\an$-admissible curve, and let $\{\phi_s\}$ denote the flow of
$\T \circ \ro \circ \s$. Then the $h^c$-transport operator
$\psi^b_a : V_{c(a)}\t \to V_{c(b)}\t$ along $\psi$ is given by
$\psi^b_a = T\phi_{b-a}$.
\end{lem}
\begin{pf}
Representing the flow of $\fpd{}{t}$ on $\R$ by $\{\l_s\}$, it
immediately follows from $\t \circ \phi_s = \l_s \circ \t$ that,
for any $v \in V\t$, the vector $T\phi_{s}(v)$ also belongs to
$V\t$. In particular, we have $T\phi_{b-a}(v) \in V_{c(b)}\t$.

Next, take $v_0 \in V_{c(a)}$ and let $X(t)$ denote the section of
$V\t$ along $c(t)$ which is uniquely determined by the conditions
$D_{\psi}(u,X)(t) = 0$ and $X(a)=v_0$. This is still equivalent to
\begin{equation}
\label{completelift} \left.\frac{d}{dt}\right|_t X(t) = \mbox{\fr
s}\,(T (\T \circ \ro) T_{c(t)}\s (X(t))).\end{equation} Since
$\mbox{\fr s}\,(T (\T \circ \ro) T_{c(t)}\s (X(t))) = (\T \circ
\ro \circ \s)^c(X(t))$, where $(\T \circ \ro \circ \s)^c$ denotes
the complete lift of the vector field $\T \circ \ro \circ \s$ to
$TM$, (\ref{completelift}) tells us that $X(t)$ is an integral
curve of $(\T \circ \ro \circ \s)^c$, passing through $v_0$. By
construction of the complete lift of a vector field, the flow of
$(\T \circ \ro \circ \s)^c$ is given by $\{T\phi_s\}$ and,
therefore, $X(t) = T\phi_{t-a}(X(a))$. The result then follows
immediately from the definition of the $h^c$-transport operator
along $\psi$.\qed
\end{pf}
Next, we consider the case where $\psi:[a,b] \to J^1\nu$ is a
piecewise $\an$-admissible curve whose projection $c = \nu_1 \circ
\psi$ onto $M$ is continuous. Recall, in particular, that
$u(t):=\nu_{1,0}(\psi(t))$ is a control (see Lemma \ref{lemma21}).
For the sequel we will need an extension of the definition of the
$\an$-derivative corresponding to $h^c$ to piecewise curves. For
that purpose, let $X: [a,b] \to V\t$ be a continuous piecewise
curve projecting onto the base curve $c(t)$ of $\psi$. Note, in
particular, that $(u,X)$ represents a piecewise section of $U
\times_M V\t$ along $c$ in $M$. From the definition of piecewise
curves it can be deduced that one can always find a sufficiently
fine subdivision $a_0=a < a_1 \ldots <a_{\ell}=b$ of the given
interval $[a,b]$ such that $\psi$ can be written as a
concatenation of $\ell$ smooth $\an$-admissible curves
$\psi_i:[a_{i-1},a_i] \to J^1\nu$ and $X$ as a concatenation of
$\ell$ smooth curves $X_i:[a_{i-1},a_i] \to V\t$. For the
piecewise $\an$-admissible curve $\psi$ we now define the
$\an$-derivative $D_{\psi}$, acting on the piecewise section
$(u,X)$, as follows:
\[ D_\psi(u,X)(t):=D_{{\psi}_i}(u,X_i)(t) \quad \mbox{for all}\;
t \in ]a_{i-1},a_i],\;i=1,\ldots,\ell,\] and \[
D_\psi(u,X)(a_0)=D_{\psi_1}(u,X_1)(a_0).\] It is easily seen in
coordinates, for instance, that the mapping $D_{\psi}(u,X)$ from
$[a,b]$ to $U \times_M V\t$ is indeed well defined. Given any $v_0
\in V_{c(a)}\t$, one can readily verify that there exists a unique
continuous piecewise curve $X(t)$ in $V\t$ such that
$D_\psi(u,X)(t) = 0$ for all $t \in [a,b]$, with $X(a)=v_0$. This
implies that one may introduce a (composite) $h^c$-transport
operator $\psi^b_a$ along the piecewise $\an$-admissible curve
$\psi$ as follows: $\psi^b_a = (\psi_\ell)^{a_\ell}_{a_{\ell-1}}
\circ \ldots \circ (\psi_1)^{a_1}_{a_0}$, where
$(\psi_i)^{a_i}_{a_{i-1}}$ represents the $h^c$-transport operator
along the smooth $\an$-admissible curve $\psi_i$, as defined in
the previous section. If, for a given $\psi$ (and the
corresponding control $u$), $X(t)$ solves the equation
$D_\psi(u,X)(t)=0$, it then follows from the definition that
$\psi^b_a(X(a)) =X(b)$.

We shall prove below that a piecewise $\an$-admissible curve
$\psi$, with $\nu_1\circ \psi = c$, can always be considered as
being induced by smooth $\an$-admissible curves
$\psi^i:I_i=[a_{i-1},a_i] \to J^1\nu$ of the form $\psi_i(t) =
j^1\s_i(c(t))$, for some local section $\s_i$ of $\nu$. Using this
property we then know from above that $(\psi_i)_{a_{i-1}}^{a_i} =
T\phi^i_{a_i-a_{i-1}}$, with $\{\phi^i_s\}$ the flow of $\T\circ
\ro \circ\s_i$. Denoting the composite flow of the ordered set
$(\T\circ \ro \circ\s_{\ell}, \ldots , \T\circ \ro \circ \s_{1})$
by $\Phi$ and using the shorthand notation introduced in Section
2, we find that the (composite) $h^c$-transport operator
$\psi^b_a$ is given by \[ \psi^b_a = T\Phi^b_a.\] Indeed, a
straightforward computation gives: \begin{eqnarray*} \psi^b_a
(X(a)) &:=& (\psi_\ell)^b_{a_{\ell-1}} \circ \ldots \circ
(\psi_1)^{a_1}_{a}(X(a)) \\
&=& T \phi^\ell_{a_\ell-a_{\ell-1}} \circ \ldots
\circ T\phi^1_{a_1-a}(X(a))\\
&=&T\Phi_a^b(X(a))=X(b). \end{eqnarray*}

In order to prove that any (piecewise) $\an$-admissible curve can
be written as a concatenation of smooth $\an$-admissible curves of
the form $j^1\s\circ c$, we shall prove that any smooth
$\an$-admissible curve $\psi$ whose image is entirely contained in
a coordinate chart, is of that form. The general result then
follows by a similar argument as the one applied in Section
\ref{sectiegeometricsetting} (when proving that the base curve of
any control is a concatenation of integral curves of vector fields
in $\cal D$). So, assume $\psi$ can be written in coordinates as
$\psi(t)=(t,x^i(t),u^a(t),u^a_t(t),u^a_i(t))$ for all $t \in
I=[a,b]$. Since $\psi$ is $\an$-admissible, we then have that
\[\dot{ u}^a(t) = u^a_t(t) + u^a_i(t) \dot x^i(t) \quad \mbox{and}
\quad \dot x^i(t) =\ro^i(t,x^i(t),u^a(t)).\] Consider now a smooth
extension $\tilde \psi(t)=(t,\tilde x^i(t),\tilde u^a(t),\tilde
u^a_t(t), \tilde u^a_i(t))$ of $\psi$, defined on an open interval
$\tilde I$ containing $I$, such that $\im \tilde\psi$ is still
contained in the same coordinate chart, with $\tilde \psi(t) =
\psi(t)$ for all $t \in I$. Next, we can construct a local section
$\s$ of $\nu$, defined on $\t^{-1}(\tilde I)$, as follows:
$\s(t,x)=(t,x,\s^a(t,x))$, with $\s^a(t,x) = \tilde u^a(t) +
\tilde u^a_i(t)(x^i - \tilde x^i(t))$. For each fixed $t \in I$ we
find that \begin{eqnarray*} \s^a(t,x^i(t)) &=&
u^a(t)\\\fpd{\s^a}{t}(t,x^i(t)) &=& \dot {u}^a(t) - u^a_i(t)
\dot { x}^i(t) =  u^a_t(t),\\
\fpd{\s^a}{x^i}(t,x^i(t)) &=& u^a_i(t),\\
\end{eqnarray*} and, hence, we have that $j^1\s(t,x(t)) = \psi(t)$ for all $t
\in I$, which is precisely what we wanted to prove.

We have seen that, given a piecewise $\an$-admissible curve $\psi$
in $J^1\nu$, with continuous piecewise base curve $c = \nu_1 \circ
\psi$ and corresponding control $u = \nu_{1,0}\circ \psi$, we can
regard the equation $D_\psi(u,X)(t)=0$ as a differential equation
for the component $X$ of the curve $(u,X)$ in $U \times_M V\t$
that is $h^c$-transported along $\psi$. Returning to the given
geometric control structure $(\t,\ro,\nu)$, we shall now explain
the role of the $h^c$-transport operator in determining the
vertical variational cone associated to a composite flow $\Phi$
and composite flow parameter $T$ induced by an ordered set of
vector fields of the form $(\T \circ \ro \circ \s)$, for some $\s
\in \G(\nu)$.

Given any control $u:[a,b] \to U$, with base curve $c = \nu \circ
u$. In Section 2 we have seen that $c$ is induced by the composite
flow $\Phi$ of an ordered set of vector fields belonging to the
family $\cal D$ given by (\ref{distr}), say $(\T \circ \ro \circ
\s_\ell,\ldots, \T \circ \ro \circ \s_1)$, where $\s_i \in
\G(\nu)$, and let the composite flow parameter be
$T=(a_\ell-a_{\ell-1},\ldots,a_1-a_0)$, with $a=a_0 < a_1 < \ldots
< b=a_\ell$. If we put $c(a_0)=m$ and $c(a_{\ell})=m'$, then the
vertical variational cone $V_{m'}R_m(\Phi,T)$ is completely
determined by the piecewise $\an$-admissible curve $\psi$ in
$J^1\nu$ that is induced by the smooth curves
$\psi_i(t)=j^1\s_i(c(t))$. Indeed, it follows from Definition
\ref{verticalcone} and from the above analysis, that any element
of $V_{m'}R_m(\Phi,T)$ can be written as a linear combination of
$h^c$-transported vertical tangent vectors along $\psi$, i.e.

\begin{eqnarray*} V_{m'}R_m(\Phi,T) &=& \left\{ \sum_i\d^i
\psi^b_{\t^i}(Y_i(c(\t^i))-\dot c(\t^i)) \ | \right.\\ && \ Y_i
\in {\cal D}, \d^i \ge 0, \t^i \in]a,b]\Bigg\}. \end{eqnarray*}

Roughly speaking, one can say that the (piecewise)
$\an$-admissible curve $\psi$ corresponding to the control $u$,
contains sufficient information regarding the sections $\s_i$ in
order to determine the vertical variational cone $V_{m'}R_m$. From
now on we shall therefore write $V_{m'}R_m(\psi)$ if we want to
emphasise that the vertical variational cone can be generated by
the $h^c$-transport operator along the (piecewise)
$\an$-admissible curve $\psi$.

For later use we will need an extension of the action of the
$\an$-derivative $D$ to `vertical' forms, belonging to the dual of
$\vvectorfields{\nu}$. Consider the fibred product bundle $U
\times_M V^*\t$ with corresponding projections $p_1^*: U \times_M
V^*\t \to U$, $p_2^*:U \times_M V^*\t \to V^*\t$, such that $\nu
\circ p^*_1 = \t^*_M \circ V^*\t$. Here $\t^*_M: V^*\t \to M$
denotes the dual bundle of $V\t \to M$. The dual module of
$\vvectorfields{\nu}$ is then given by the set
\[ \voneforms{\nu}=\{\eta: U \to V^*\t\,|\,\t^*_M(\eta(u)=\nu(u)\;
\mbox{for all}\; u \in U\}.\] Obviously, we have
$\voneforms{\nu}\cong \G(p^*_1)$.  Given $\eta \in
\voneforms{\nu}$ and $Z \in \vvectorfields{\nu}$, the natural
pairing $\vvectorfields{\nu}$, $\langle \eta, Z \rangle$ defines a
function on $U$. In particular, for $j^1_m\s \in J^1\nu$, with
$\s(m)=u$, we note that $p_2^*(\eta(u))$ and $p_2(Z(u))$ belong to
the dual linear spaces $V_m^*\t$ and $V_m\t$, respectively. By
requiring that for any fixed $\eta \in \voneforms{\nu}$, the
relation
\begin{equation}\label{deriv}
\langle D_{j^1_m\s} \eta , Z(u) \rangle = \an(j^1_m\s) (\langle
\eta,Z\rangle) - \langle \eta(u), D_{j^1_m\s}Z\rangle,
\end{equation}
should hold for all $Z \in \vvectorfields{\nu}$, the element
$D_{j^1_m\s}\eta \in (p^*_1)^{-1}(u) (\cong V^*_m\t)$ is uniquely
determined.

Consider a piecewise $\an$-admissible curve $\psi:[a,b]\to J^1\nu$
with continuous piecewise projection $c=\nu_1 \circ \psi$ on $M$
and corresponding control $u = \nu_{1,0} \circ \psi: [a,b] \to U$.
Take a continuous piecewise section $\overline{\eta}(t)$ of
$V^*\t$ along $c(t)$ such that $(u,\overline{\eta})(t)$ defines a
section of $p_1^*$ along the curve $u(t)$. We then have the
following property.

\begin{lem}$D_\psi (u,\overline{\eta})(t) =0$ iff
$\overline{\eta}(t)=((\psi^t_a)^{-1})^*(\overline{\eta}(a))$ for
all $t \in [a,b]$.\end{lem}
\begin{pf} Fix some $t_0 \in I$ and take an arbitrary $X_0 \in
V_{c(t_0)}\t$. Using the $h^c$-transport operator along $\psi$, we
can then construct a continuous piecewise section $X(t)$ of $V\t$
along $c(t)$ by $X(t) = \psi^{t}_a((\psi^{t_0}_a)^{-1}(X_0))$.
Note that $X(t_0)=X_0$. Then, with (\ref{deriv}) we obtain
\begin{eqnarray*}\langle D_\psi (u,\overline{\eta})(t_0), (u(t_0),X(t_0))
\rangle &=& \left.\frac{d}{dt}\right|_{t_0} \langle
\overline{\eta}(t),X(t)\rangle
\\&& - \langle
(u(t_0),\overline{\eta}(t_0)),D_\psi(u,X)(t_0)\rangle.\end{eqnarray*}
Now it follows from the definitions that both terms on the
right-hand vanish separately if we take
$\overline{\eta}(t)=((\psi^t_a)^{-1})^*(\eta(a))$. Indeed, with
this choice we have $\langle \overline{\eta}(t),X(t)\rangle \equiv
\langle \eta(a),(\psi^{t_0}_a)^{-1}(X_0)\rangle =\,
\mbox{const.}$, and the equation $D_\psi(u,X)(t_0)=0$ holds in
view of the definition of $X(t)$. The remainder of the proof then
follows from the uniqueness of solutions of a system of ordinary
differential equations with given initial conditions.\qed
\end{pf}
The $\an$-derivative $D$ will play a crucial role in the proof of
the Maximum Principle in the next section. In the following remark
we briefly explain how some of the basic ideas in the treatment of
the Maximum Principle in \cite{Pont} can be related to our work.
\begin{rem} \textnormal{
The discussion of the Maximum Principle can be developed for
controls that verify the weaker assumption of being measurable and
bounded, instead of (piecewise) smooth (see, for instance, L.S.
Pontryagin \etal \ \cite{Pont}). Using local coordinate
expressions, we will roughly sketch how the smoothness conditions
we have imposed on controls can also be relaxed within our
framework. The local expressions for the equation $D_\psi (u,X)(t)
= 0$ reads
\[ \dot X^k(t)=
\left(\fpd{\ro^k}{x^i}(t,c^j(t),u^a(t)) +
u^b_i(t)\fpd{\ro^k}{u^b}(t,c^j(t),u^a(t))\right)X^i(t).
\]
The condition that the functions $u^a(t)$ and $u^a_i(t)$ be
measurable and bounded, suffices to obtain a solution of this
equation and, subsequently, to introduce a suitable notion of
transport operator.
This observation can be translated into our geometric framework as
follows. Consider the set $V^1\nu := \cup_{m \in M}\{ T_m\s_{|V\t}
: V_m\t \to T_{\s(m)}U \ | \  \s \in \G(\nu)\}$. It can be proven
by standard arguments that $V^1\nu$ is an affine bundle over $U$,
with coordinates $(t,x^i,u^a,u^a_i)$ (see, for instance,
\cite{saunders}). Note that there exists a natural projection
$\mu:J^1\nu \to V^1\nu$, locally expressed by
$(t,x^i,u^a,u^a_i,u^a_t) \mapsto (t,x^i,u^a,u^a_i)$. From the fact
that the coefficients $\G^k_i$ of $h^c$ do not depend on the
$u^a_t$ (see the proof of Proposition 6.3) it easily follows that
the $\an$-derivative $D_\psi$ only depends on $\mu \circ \psi$.
Now, since $\psi$ was assumed to be $\an$-admissible, i.e.
$\an(\psi)=\dot u$, the smoothness condition on $u$ could not be
relaxed. However, the curve $\tilde \psi = \mu \circ \psi$ does
not have to satisfy this condition, implying that the smoothness
condition can be relaxed without losing the notion of derivative
acting on sections of $V\t$ along $c$. We can therefore conclude
that, in order to define a vertical cone of variations associated
with a measurable and bounded control $u$, we must fix a curve
$\tilde \psi$ in $V^1\nu$.  If one works in a coordinate chart, a
natural choice of $\tilde \psi$ is the curve $\tilde
\psi(t)=(t,c^i(t),u^a(t),u^a_i(t))$ with $u^a_i(t)=0$. The
equations of the derivative associated with $\tilde \psi$ then
reduce to $\dot X^k(t) = \fpd{\ro^k}{x^i}(t,c^j(t),u^a(t))X^i(t)$.
These equations are precisely the ``variational equations''
introduced in \cite[p79]{Pont}. By fixing the coordinate chart,
one can fix the section $\s^a(t,x) = u^a(t)$ and the curve $\tilde
\psi(t)=(t,c^j(t),u^a(t),0)$, implying that, respectively a fixed
vertical cone of variations and a fixed derivative associated with
$\tilde \psi$ can be defined. This essentially establishes the
link between our approach and the one followed by L.S. Pontryagin
\etal.}
\end{rem}

\section{The Maximum Principle and extremal controls}
We will now derive the Maximum Principle by combining the tools
developed in Section \ref{sectiecontrolafgeleide} and the
necessary conditions for optimal controls derived in Section
\ref{costoptimal}.

Let $(\t,\nu,\ro,L)$ denote an arbitrary geometric optimal control
structure, with extended geometric control structure $(\ovl
\t,\ovl \nu,\ovl \ro)$. In view of the structure of the bundle
$\ovl \t: \ovl M (= M \times \R) \to \R,\; (m,J) \mapsto \t(m)$,
it is easily seen that the bundle of vertical tangent vectors
$V\ovl \t$ is isomorphic to $V\t\times \R^2$. Similarly, the
bundle $V^*\ovl \t$ can be identified with $V^*\t \times \R^2$. In
particular, given a point $(m,J) \in \ovl M (= M \times\R)$, a
co-vector $\ovl \eta \in V^*_{(m,J)}\ovl \t$ can always be
represented by a pair $(\eta_m,\eta_J)$ for some $\eta_m \in
V^*_m\t$ and $\eta_J \in \R$.

Before proceeding, we still have to introduce a few additional
concepts. First, we recall that the dual of a convex cone $C$ in a
vector space ${\cal V}$ is defined by the set $C^*=\{\a \in {\cal
V}^* \ | \ \langle \a,v\rangle \le 0 ,\ \forall v \in C\}$. A
general result that will be used later on, tells that
$C^*=(\mbox{cl}(C))^*$ and $(C^*)^*=\mbox{cl}(C)$, where
$\mbox{cl}$ denotes the closure of $C$ in ${\cal V}$ (see e.g.
\cite{kothe} for a proof). Finally, for any $v \in {\cal V}$, the
half-ray through $0$ en $v$, i.e. $\{w \ |\ w = r v, \ \forall r
\ge 0\}$, will be called the `cone generated by $v$', and denoted
$C(v)$.

Another concept that we will need, is that of a `multiplier of a
control'. For that purpose, we first construct a $1$-parameter
family of closed two-forms on $U\times_M V^*\t$. Let $\tilde \w$
be the closed two-form on the fibred product $U\times_M T^*M$,
obtained by pulling back the canonical symplectic form on $T^*M$
by the projection $U\times_M T^*M\to T^*M$. Next, for any real
number $\l$ we can define a section $H_\l$ of the fibration
$U\times_MT^*M \to U \times_M V^*\t$ in the following way. Take $u
\in U_m, \eta \in V^*_m\t$ and put $H_{\l}(u,\eta) =(u, \a)$,
where $\a \in T^*_mM$ is uniquely determined by the conditions
$\langle \a, \T(\ro(u))\rangle + \l L(u) = 0$ and $\a$ projects
onto $\eta$. The mapping $H_{\l}$ is smooth, as can be easily seen
from the following coordinate expression: putting $u=(t,x^i,u^a)$
and $\eta = p_i\,dx^i_{|m}$, a straightforward computation gives

\[H_{\l}(t,x^i,u^a,p_i) = \left(t,x^i,u^a,-\ro^i(t,x^i,u^a)p_i
-\l L(t,x^i,u^a), p_i \right).
\]
We can now use $H_{\l}$ to pull-back the closed two-form $\tilde
\w$ to a closed two form on $U\times_M V^*\t$, which will be
denoted by $\w_{\l}$. Herewith, we can now introduce the following
definition of a multiplier.
\begin{defn}\label{multiplier}
Given a control $u: [a,b] \to U$, a pair $(\eta,\l)$ consisting of
a continuous piecewise section $\eta$ of $V^*\t$ along $c = \nu
\circ u$ and a real number $\l$, is called a {\rm multiplier} of
$u$ if the following conditions are satisfied:
\begin{enumerate}
\item \label{max1}$i_{(\dot u(t),\dot \eta(t))} \w_{\l} =0$ on
every smooth part of the curve $(u(t),\eta(t))$, \item
\label{max2}given any $t_0 \in [a,b]$, and putting
$H_{\l}(u(t_0),\eta(t_0))= (u(t_0),\a_0)$, the function
$u'\mapsto\langle \a_0, \T(\ro(u')) \rangle + \l L(u')$, defined
on $\nu^{-1}(c(t_0))$, attains a global maximum for $u'=u(t_0)$,
\item \label{max3} $(\eta(t),\l) \neq (0,0)$ for all $t \in [a,b]$.
\end{enumerate} \end{defn}

Returning to the geometric optimal control structure
$(\t,\nu,\ro,L)$, let $\ovl u(t) = (u(t),J_0+{\cal J}_u^{(a,t)})$
represent a control in the extended geometric control setting,
defined on an interval $[a,b]$. As before, $c$ will denote the
base curve of $u$ in $M$ (cf. Section \ref{costoptimal}), and we
put $c(a)=m, c(b)=n$. The bundle map (\ref{lambda}) associated to
the extended geometric control structure will be written as $\ovl
\an$. Given an arbitrary piecewise $\ovl \an$-admissible curve
$\ovl \psi$ in $J^1\ovl \nu$ projecting onto $\ovl u$, we will
prove in the following theorem that the dual of the vertical
variational cone $V_{(n,J_0+J(u))}R_{(m,J_0)}(\ovl \psi)$ only
depends on $u$.

\begin{thm}
Let $\ovl \eta_0 =(\eta_0,\l_0) \in V^*_{(n,J_0 + {\cal J}(u))}
\ovl \t$, with $\ovl \eta_0 \neq 0$. Then we have that $\ovl
\eta_0 \in \left(V_{(n,J_0+{\cal J}(u))}R_{(m,J_0)}(\ovl
\psi)\right)^*$ if and only if there exists a section $\eta$ of
$V^*\t$ along $c$, with $\eta(b) = \eta_0$, such that the pair
$(\eta,\l_0)$ is a multiplier of $u$.
\end{thm}
\begin{pf}
We prove that any $\ovl \eta_0 \neq 0$ in the dual of the vertical
variational cone determines a multiplier for $u$. The converse
property will then simply follow by reversing the arguments.

Let $\ovl \eta(t)$ denote the unique continuous piecewise curve in
$V^*\ovl\t$ satisfying the equation $\ovl D_{\ovl \psi} (\ovl
u,\ovl \eta)(t) =0$, with $\ovl \eta(b) = \ovl \eta_0$. This
implies that $\ovl \eta(t) = (\psi^b_t)^{*} (\ovl \eta_0)$. We can
write $\ovl \eta(t)$ as $\ovl \eta(t) = (\eta(t),(J_0 + {\cal
J}_{u}^{(a,t)},\eta_J(t)))$, where $\eta(t)$, resp.\ $(J_0 + {\cal
J}_{u}^{(a,t)},\eta_J(t))$ are curves in $V^*\t$, resp.\ $\R^2$,
such that $\eta(b)=\eta_0$ and $\eta_J(b)=\l_0$. We will now prove
that $(\eta,\l_0)$ is a multiplier of $u$.

First of all, it is easily seen that condition (\ref{max3}) of
Definition \ref{multiplier} holds. In order to prove that (1) and
(2) of the definition hold, take an arbitrary $t_0 \in ]a,b]$ and
$u' \in U_{c(t_0)}$ arbitrary. Then, we find that
\[
\ovl \psi^b_{t_0}\left( \T(\ovl \ro(u',J_0+{\cal J}_u^{(a,t_0)}))
- \T(\ovl \ro(\ovl u(t_0))) \right) \in V_{(n,J_0+{\cal
J}(u))}R_{(m,J_0)}(\ovl \psi).
\]
By contracting this tangent vector with $\ovl \eta_0 \in
\left(V_{(n,J_0+{\cal J}(u))}R_{(m,J_0)}(\ovl \psi)\right)^*$, and
taking into account the definition of the dual of a cone, we
obtain the following inequality:
\begin{equation}\label{gm}
\langle \eta(t_0) , \T (\ro(u')) -\T( \ro(u(t_0))) \rangle +
\eta_J(t_0) (L(u') - L(u(t_0))) \le 0.
\end{equation}
This holds for any $t_0 \in ]a,b]$ and any $u' \in U_{c(t_0)}$.
Note that this inequality is also valid for $t_0=a$. It suffices
to consider a local trivialization of $U$ and to interpret the
left-hand side of the above inequality as a function of $t_0$,
which is clearly continuous in a neighborhood of $a$. In
particular, we deduce from the above that the function
\[
u' \mapsto \langle \eta(t_0) , \T (\ro(u')) -\T( \ro(u(t_0)))
\rangle + \eta_J(t_0) (L(u') - L(u(t_0))),
\]
defined on $U_{c(t_0)}$, admits a global maximum at $u' = u(t_0)$.
In local coordinates this means, in particular, that we have:
\begin{equation}\label{extremum}
\eta_i(t_0)\fpd{\ro^i}{u^a}(u(t_0)) + \eta_J(t_0)
\fpd{L}{u^a}(u(t_0)) = 0,
\end{equation}
and this holds for all $t_0 \in [a,b]$. These relations are used
in the following to prove that the function $\eta_J(t)$ is
constant and that $\eta(t)$ satisfies condition (\ref{max1}) of
Definition \ref{multiplier}. The coefficients of the linear $\ovl
\an$-lift ${\ovl h}^c$ are related to the coefficients of $h^c$ in
the following way (using a slight abuse of notation):
\[ \begin{array}{ll} \ovl \G^i_j = \G^i_j, & \displaystyle\ovl \G^i_J =
\fpd{\ro^i}{u^a} u^a_J, \vspace{.2cm}\\ \displaystyle \ovl \G^J_i
= \fpd{L}{x^i} + \fpd{L}{u^a}u^a_i, & \displaystyle \ovl \G^J_J =
\fpd{L}{u^a}u^a_J.
\end{array}\]
Herewith, the differential equations for $\eta(t)$ and $\eta_J(t)$
become, on every smooth part of $\ovl \eta$:
\begin{eqnarray*} \dot \eta_J(t) &=& - \ovl \G_J^i \eta_i(t) - \ovl \G^J_J\eta_J(t),\\
\dot \eta_i(t) &=& - \G^j_i \eta_ j(t) - \ovl
\G^J_i\eta_J(t).\end{eqnarray*} Taking into account the relations
(\ref{extremum}), which hold for all values of $t_0 \in [a,b]$, it
is easily seen that $\dot{\eta}_J(t)= 0$ and, hence, $\eta_J$ is a
constant function, with $\eta_J (t) \equiv \l_0$. Moreover, the
functions $\eta_i(t)$ satisfy:
\[
\dot \eta_i(t) = - \fpd{\ro^j}{x^i} \eta_j(t) -\l_0 \fpd{L}{x^i}.
\]
Putting, in local coordinates, $h_{\l_0}(u,\eta) = \ro^i(u)\eta_i
+ \l_0 L(u)$, the two-form $\w_{\l_0}$ reads: $\w_{\l_0} =
-dh_{\l_0} \wedge dt + dp_i \wedge dx^i$. After some tedious, but
straightforward calculations it follows that the condition
$i_{(\dot u(t), \dot \eta(t))} \w_{\l_0} = 0$ is equivalently to
\begin{eqnarray*}
&&\dot c^i(t) = \fpd{h_{\l_0}}{p_i}(u(t),\eta(t))=\ro^i(u(t)),\\
&&0= \fpd{h_{\l_0}}{u^a}(u(t),\eta(t))=
\fpd{\ro^i}{u^a}(u(t))\eta_i(t) + {\l_0} \fpd{L}{u^a}(u(t)),\\
&& \dot \eta_i(t) = - \fpd{h_{\l_0}}{x^i}(u(t),\eta(t))=
-\fpd{\ro^j}{x^i}(u(t))\eta_j(t) - {\l_0} \fpd{L}{x^i}(u(t)),\\
&& \left.\frac{d}{dt}\right|_t\left( h_{\l_0}(u(t),\eta(t))
\right) = \fpd{h_{\l_0}}{t}(u(t),\eta(t)),
\end{eqnarray*}
and it is easily seen that the curve $\eta(t)$ defined above,
satisfies these equations. This shows that condition (3) of
Definition \ref{multiplier} is satisfied. It finally remains to
prove that also the second condition for a multiplier holds.

Consider the section $H_{\l_0}: U \times_M V^*\t \to U\times_M
T^*M$ and let us write for any $t_0 \in [a,b]$,
$H_{\l_0}(u(t_0),\eta(t_0))=(u(t_0),\a(t_0))$. Substituting this
into (\ref{gm}), and recalling that $\eta_J(t)\equiv \l_0$, we
obtain:
\[
\langle \a(t_0), \T(\ro(u')) \rangle  + \l_0 L(u')\left(=\langle
\a(t_0) , \T(\ro(u(t_0))) \rangle + \l L(u(t_0))\right) \le 0 ,
\]
proving that (\ref{max2}) is satisfied. This completes the proof
that $(\eta(t),\eta_J =\l_0)$ is indeed a multiplier. \qed
\end{pf}
As a consequence of the above theorem, the dual of the vertical
variational cone, in the extended setting, only depends on the
control $u$ and, hence, this is also true for the closure of this
cone. Moreover, as an interesting side result we obtain that the
closure of the vertical variational cone $V_nR_m$ also depends on
$u$ only. Indeed, using the same techniques as in the above
theorem it is easily seen that every multiplier with $\l=0$,
determines an element of the dual cone of $V_nR_m$, and vice
versa. To simplify the notations we put $\ovl m=(m,J_0)$ and $\ovl
n=(n,J_0+ {\cal J}(u))$. Recall Corollary \ref{maximum1}, which is
reformulated in the following way and leads us to a more familiar
version of the maximum principle.
\begin{cor} \label{maximum4} Assume that $m\stackrel{u}{\to} n$ and that
$u$ is optimal. Then there exists a multiplier $(\eta, \l)$ with
$\l \le 0$.
\end{cor}
The following definitions are well known from the literature.
\begin{defn} \label{extremal}
A control $u$, with $m \stackrel{u}{\to} n$ is called an extremal
if there exists a multiplier $(\eta(t),\l)$ for which $\l \le 0$.
An extremal is called normal, resp. abnormal, if there exists a
multiplier $(\eta(t),\l)$ for which $\l < 0$, resp. $\l =0$.
\end{defn}
An extremal is thus equivalently defined as a control for which
the closed cone $\mbox{cl}(V_{\ovl n}R_{\ovl m})$ does not contain
$-\partial / \partial J$ in its interior. Note that an extremal
can be simultaneously abnormal and normal. We say that an extremal
is {\em strictly abnormal} if it is abnormal but not normal. The
following proposition gives necessary and sufficient conditions
for a control to be an abnormal extremal or a strictly abnormal
extremal.
\begin{prop}
A control is an abnormal extremal iff $\mbox{\rm cl}(V_{n}R_{ m})
\neq V_n\t$. A control is a strictly abnormal extremal iff
$-\fpd{}{J}$ is in the border of $\mbox{\rm cl}(V_{\ovl n}R_{\ovl
m})$.
\end{prop}
\begin{pf}
The first statement follows from the fact that every element in
the dual cone $(\mbox{cl}(V_{n}R_{m}))^*$ corresponds to a
multiplier with $\l=0$ (see above).

An extremal is strictly abnormal iff every element $\ovl \eta_0$
in $(\mbox{cl}(V_{\ovl n}R_{\ovl m}))^*$ satisfies $(\ovl
\eta_0)_J \ge 0$ (by definition). Using the definition of the dual
cone and the fact that $C^{**} = \mbox{cl}(C)$ for an arbitrary
convex cone $C$, we obtain that $-\fpd{}{J}$ is contained in
$\mbox{cl}(V_{\ovl n}R_{\ovl m})$. On the other hand, since $u$ is
an extremal we know that $- \fpd{}{J}$ is not contained in the
interior of the cone $\mbox{cl}(V_{\ovl n}R_{\ovl m})$. \qed
\end{pf}

It should be noted that the condition $V_nR_m\neq V_n\t$ does {\em
not} depend on the cost function $L$. This justifies the notion of
an abnormal extremal: $u$ satisfies the necessary conditions for
being a optimal control with respect to the cost $L$, however
these conditions do not depend on $L$. The above result can be
intuitively interpreted as follows: a control $u$ is an abnormal
extremal iff the family of vector fields ${\cal D}$ does not
supply enough ``vertical'' variations to the control $u$. In the
case of strictly abnormal extremals the maximum principle fails in
the sense that Corollary \ref{approx2} only gives information on
those vectors lying in the interior of a variational cone, and not
on those belonging to the boundary.

\begin{acknowledgements}
This work has been supported by a grant from the ``Bijzonder
Onderzoeksfonds'' of Ghent University. I am indebted to F.
Cantrijn for the many discussions and the careful reading of this
paper and to J. Cort\'es and A. Ibort for many useful suggestions.
\end{acknowledgements}

\end{article}
\end{document}